\documentclass[11pt,twoside,a4paper]{article}

\usepackage{fullpage}
\usepackage{enumerate}
\usepackage{subfigure}
\usepackage{xcolor}

\usepackage[utf8]{inputenc}
\usepackage[T1]{fontenc}
\usepackage[english]{babel}
\usepackage{enumitem}
\usepackage{amsmath,amsfonts,amssymb,graphicx,bbm}
\usepackage{stackrel}

\usepackage{amsthm}
\usepackage{constants}

\theoremstyle{plain}
\newtheorem{theorem}{Theorem}

\newtheorem{proposition}{Proposition}[section]
\newtheorem{lemma}[proposition]{Lemma}
\newtheorem{corollary}[proposition]{Corollary}

\newtheorem{conjecture}{Conjecture}

\theoremstyle{definition}

\newtheorem{remark}[proposition]{Remark}

\usepackage{mathtools}
\DeclarePairedDelimiter\abs{\lvert}{\rvert}

\newcommand{\mbb}[1]{\mathbb{#1}}
\newcommand{\mbf}[1]{\mathbf{#1}}
\newcommand{\mcal}[1]{\mathcal{#1}}
\newcommand{\msf}[1]{\mathsf{#1}}

\newcommand{\mrm}[1]{\mathrm{#1}}

\newcommand{\bbm}[1]{\mathbbm{#1}}

\newcommand{\Z}{\mathbb{Z}}
\newcommand{\R}{\mathbb{R}}

\newcommand{\Prob}[1]{\mathbf{P}\left[#1\right]}

\newcommand{\at}{\mathtt{AT}}
\newcommand{\hf}{\mathtt{HF}}
\newcommand{\gat}{\mathtt{GAT}}
\newcommand{\atrc}{\mathtt{ATRC}}
\newcommand{\eightv}{\mathtt{EightV}}
\newcommand{\sixv}{\mathtt{SixV}}

\newcommand{\spin}{s}

\newcommand{\ind}[1]{\mathbbm{1}_{\{#1\}}}

\newcommand{\norm}[1]{\left\Vert #1 \right\Vert}

\usepackage{verbatim}

\usepackage{tabularx}

\usepackage{hyperref}
\hypersetup{
    colorlinks=true,
    linkcolor=blue,
    citecolor=red,
    urlcolor=blue,
    pdfborder={0 0 0}
}

\title{On antiferromagnetic regimes in the Ashkin--Teller model}

\date{\today}

\author{Moritz Dober\thanks{University of Vienna, \url{moritz.dober@univie.ac.at}}}

\begin{document}

\maketitle

\begin{abstract}
The Ashkin--Teller model can be represented by a pair $(\tau,\tau')$ of Ising spin configurations with coupling constants $J$ and $J'$ for each, and $U$ for their product.
We study this representation on the integer lattice $\mathbb{Z}^d$ for $d\geq 2$.
We confirm the presence of a partial antiferromagnetic phase in the isotropic case ($J=J'$) when $-U>0$ is sufficiently large and $J=J'>0$ is sufficiently small, by means of a graphical representation.
In this phase, $\tau$ is disordered, admitting exponential decay of correlations, while the product $\tau\tau'$ is antiferromagnetically ordered, which is to say that correlations are bounded away from zero but alternate in sign.
No correlation inequalities are available in this part of the phase diagram.
In the planar case $d=2$, we construct a coupling with the six-vertex model and show, in analogy to the first result, that the corresponding height function is localised, although with antiferromagnetically ordered heights on one class of vertices of the graph.

We then return to $d\geq 2$ and consider a part of the phase diagram where $U<0$ but where correlation inequalities still apply. Using the OSSS inequality, we proceed to establish a subcritical sharpness statement along suitable curves covering this part, circumventing the difficulty of the lack of general monotonicity properties in the parameters. We then address the isotropic case and provide indications of monotonicity.
\end{abstract}

\section{Introduction}
\subsection{The Ashkin--Teller and related models}
The \emph{Ashkin--Teller} (AT) model~\cite{AshTel43} may be represented~\cite{Fan72} by a pair $(\tau,\tau')$ of interacting Ising spin configurations, which are random assignments of $\pm1$ \emph{spins} to the vertices of some graph  $\mbb{G}=(\mbb{V}_\mbb{G},\mbb{E}_\mbb{G})$.
In its general form, the \emph{spin representation} involves three \emph{coupling constants} $J,J',U\in\R$, where $J$ and $J'$ describe the strength of interaction within each of the Ising configurations $\tau$ and $\tau'$, while $U$ describes the strength of interaction for their product $\tau\tau'$ and hence the relation between $\tau$ and $\tau'$. Its formal \emph{Hamiltonian} is given by 
\begin{equation*}
-\sum_{\{x,y\}\in\mbb{E}_{\mbb{G}}}J\tau_x\tau_y+J'\tau'_x\tau'_y+U\tau_x\tau'_x\tau_y\tau'_y.
\end{equation*}
When $J$ and $J'$ coincide, the model is referred to as \emph{isotropic}, and \emph{anisotropic} otherwise.
In the case $U=0$, it reduces to two independent Ising models. For $J=J'=U$, the 4-state Potts model is recovered.
In the \emph{ferromagnetic} case, where all coupling constants are non-negative, it was conjectured in~\cite{Weg72} and later on in~\cite{WuLin74} that the AT model on the square lattice $\Z^2$ generally undergoes two \emph{phase transitions}. This was verified when $U>J+J'$~\cite{Pfi82} and when $U>J=J'$~\cite{AouDobGla24}. The phase diagram of the isotropic model on the square lattice was predicted in~\cite{Kno75,DitBanGreKad80}, see also~\cite[Chapter 12.9]{Bax89} and~\cite{HuaDenJacSal13,GlaPel23}. Its ferromagnetic part was established in~\cite{AouDobGla24}.
The presence of negative coupling constants results in preferential disagreement along the edges, which is referred to as \emph{antiferromagnetic} interaction. This further leads to a loss of correlation inequalities and monotonicity properties in the parameters, thereby significantly complicating the analysis of the model.

In case the underlying graph is \emph{bipartite}, such as $\Z^d$, performing spin flips for $\tau,\,\tau'$ or both on one class of the vertices leads to a transformation of the AT measures by reversing the signs of any two of the coupling constants, while leaving the third one unchanged. 
Consequently, the main features of the model when $J$ or $J'$ are negative can be derived from those when both $J$ and $J'$ are positive.
The antiferromagnetic regime to be studied is therefore ${J,J'>0>U}$, which is the focus of this article.
In this regime, the product $\tau\tau'$ is favoured to disagree along the edges, while both $\tau$ and $\tau'$ are favoured to agree separately. 
When $J,J'<\abs{U}$, the ferromagnetic effect is dominated by the antiferromagnetic one, and one may therefore expect to observe antiferromagnetic behaviour. 
In the current article, the model is shown to exhibit an antiferromagnetic phase for $\tau\tau'$ when $J=J'$ is small and $\abs{U}$ is large, while also displaying ferromagnetic behaviour through a \emph{sharp order-disorder} phase transition for $\tau$ when $J,J'>\abs{U}$.

\paragraph{Planarity and relation to the eight-vertex model.}
If the underlying graph is \emph{planar}, duality transformations~\cite{MitSte71,Fan72c} applied to two of the AT spin configurations $\tau,\tau'$ and the product $\tau\tau'$ relate the AT models on the graph and on its \emph{dual} to one another. 
Furthermore, on the square lattice $\Z^2$, by keeping one of the AT spin configurations fixed and applying a duality transformation to one of the other, a connection~\cite{Fan72,Fan72b,Weg72,GlaPel23} is established with the spin representation~\cite{Wu71,KadWeg71,GlaMan21,Lis22,GlaPel23} of the \emph{eight-vertex} model~\cite{Sut70,FanWu70} on the \emph{medial graph}. Unless the initial AT model is self-dual, the weights of the corresponding eight-vertex model differ on two sublattices, in which case it is called \emph{staggered}~\cite{HsuLinWu75}.
In special cases when two of the AT coupling constants coincide, a certain weight vanishes in the related eight-vertex model and it reduces to the (staggered) \emph{six-vertex} model~\cite{Pau35,Rys63,WuLin75}. The latter admits a \emph{height function} representation, which has been studied in the self-dual isotropic case~\cite{DumKarManOul20,Lis21b,GlaPel23,GlaLam23}. In this particular case, the Baxter--Kelland--Wu coupling~\cite{BaxKelWu76} further relates the corresponding non-staggered six-vertex model to FK-percolation~\cite{ForKas72}, see~\cite{GlaPel23,AouDobGla24}.

\subsection{Definition of the models}\label{sec:def_models}
In this section, we introduce the Ashkin--Teller and the six-vertex height function models in a way that allows us to present our main results in the following section. See Section~\ref{sec:setting} for precise definitions regarding graphs and configurations and Section~\ref{sec:at_8v_coupling} regarding duality. 

\paragraph{The Ashkin--Teller model.}
We consider the spin representation~\cite{Fan72} of the Ashkin--Teller model~\cite{AshTel43} on the integer lattice $\Z^d$ whose edge-set we denote by $\mbb{E}_{\Z^d}$.
Define the spin configuration space $\Sigma_{\Z^d}:=\{-1,1\}^{\Z^d}$. Fix coupling constants $J,J',U\in\R$. Let $\Lambda\subset\Z^d$ be finite, and let $\eta,\eta'\in\Sigma_{\Z^d}$. Define $\Sigma_\Lambda^\eta$ as the set of all $\tau\in\Sigma_{\Z^d}$ that coincide with $\eta$ on $\Z^d\setminus\Lambda$.
The Ashkin--Teller (AT) measure on $\Lambda$ with coupling constants $J,J',U$ and \emph{boundary condition} $(\eta,\eta')$ is the probability measure on $\Sigma_{\Z^d}\times\Sigma_{\Z^d}$ given by 
\begin{equation}\label{eq:def_at}
\begin{aligned}
\at_{\Lambda,J,J',U}^{\eta,\eta'}[\tau,\tau']:=\frac{1}{Z}\exp\bigg(\sum_{\substack{\{x,y\}\in\mbb{E}_{\Z^d}:\\ \{x,y\}\cap\Lambda\neq\varnothing}}J\tau_x\tau_y+J'\tau'_x\tau'_y+U\tau_x\tau'_x\tau_y\tau'_y\bigg)\ind{(\tau,\tau')\in\Sigma_\Lambda^\eta\times\Sigma_\Lambda^{\eta'}},
\end{aligned}
\end{equation}
where the normalising constant $Z=Z(\Lambda,J,J',U,\eta,\eta')$ is called the \emph{partition function}. We write $\langle\,\cdot\,\rangle_{\Lambda,J,J',U}^{\eta,\eta'}$ for the expectation operator with respect to the above measure.
When $\eta$ or $\eta'$ are constant $+1$ or $-1$, we simply write $+$ or $-$ in the superscript, respectively. In the isotropic case, when $J=J'$, we omit $J'$ from the subscripts and simply write $\at_{\Lambda,J,U}^{\eta,\eta'}$ and $\langle\,\cdot\,\rangle_{\Lambda,J,U}^{\eta,\eta'}$.

We will also consider \emph{free} boundary conditions and define $\Sigma_\Lambda^{\mrm{f}}$ as the set of all $\tau\in\Sigma_{\Z^d}$ that are constant $1$ on $\Z^d\setminus\overline{\Lambda}$, where $\overline{\Lambda}$ is the set of vertices in $\Z^d$ that belong to $\Lambda$ or that are adjacent to vertices in $\Lambda$. Notice that the values of $\tau,\tau'$ on $\Z^d\setminus\overline{\Lambda}$ have no influence on the exponential in~\eqref{eq:def_at}.
We then also allow the boundary conditions $\eta=\mrm{f}$ and $\eta'=\mrm{f}$ in~\eqref{eq:def_at}.
It should be noted that this definition of free boundary conditions is slightly different from the standard one in the Ising model, which we comment on in Remark~\ref{rem:free_bc}.

\paragraph{The (staggered) six-vertex height function.} Consider the dual square lattice with vertex-set $(\Z^2)^*:=(\tfrac{1}{2},\tfrac{1}{2})+\Z^2$ and edges between nearest neighbours.
Regarded as line-segments between their endvertices, each edge $e$ of $\Z^2$ intersects a unique edge $e^*$ of $(\Z^2)^*$, and we call the set $e\cup e^*$ a \emph{quad}.
A (six-vertex) height function is a function $h:\Z^2\cup (\Z^2)^*\to\Z$ satisfying
\begin{enumerate}[label=(\roman*)]
\setlength\itemsep{0.1em}
\item for $x\in\Z^2$ and $x'\in (\Z^2)^*$ belonging to the same quad, $\abs{h_x-h_{x'}}=1$,
\item $h$ takes even values on $\Z^2$: for all $x\in\Z^2$, $h_x\in 2\Z$.
\end{enumerate}
We write $\Omega_{\mrm{hf}}$ for the set of all height functions.
Denote by $h^\bullet$ and $h^\circ$ the restrictions of $h$ to $\Z^2$ and $(\Z^2)^*$, respectively. For a height function $h$ and a set $E$ of edges of $\Z^2$, the sets of disagreement edges of $h^\bullet$ and $h^\circ$ in $E$ are defined respectively as
\begin{equation}\label{eq:hf_disag_edges}
\begin{aligned}
E_{h^\bullet}&:=\{e\in E:e=\{x,y\}\text{ and }h^\bullet_x\neq h^\bullet_y\},\\
E_{h^\circ}&:=\{e\in E:e^*=\{x',y'\}\text{ and }h^\circ_{x'}\neq h^\circ_{y'}\}.
\end{aligned}
\end{equation}
Observe that condition (i) implies that $E_{h^\bullet}$ and $E_{h^\circ}$ are disjoint. 
Fix $n\geq 1$, and define the boxes $\msf{B}_n:=\{-n,\dots,n\}^2\subset\Z^2$ and $\msf{B}_n^*:=\{-n+\tfrac{1}{2},\dots,n-\tfrac{1}{2}\}^2\subset(\Z^2)^*$.
Given a height~function~$t$ and weights $\mbf{a},\mbf{b},\mbf{c}\in (0,\infty)$, define the corresponding (staggered) six-vertex height function probability measure on $\Omega_\mrm{hf}$ by

\begin{equation}\label{eq:hf_meas_def}
\hf_{\msf{B}_n,\mbf{a},\mbf{b},\mbf{c}}^{t}[h]:=\frac{1}{Z}\,\mbf{a}^{\abs{E_{h^\circ}}}\,\mbf{b}^{\abs{E_{h^\bullet}}}\,\mbf{c}^{\abs{E\setminus (E_{h^\bullet}\cup E_{h^\circ})}}\,\ind{\forall x\in(\msf{B}_n\cup\mathsf{B}_n^*)^\mrm{c}:\,h_x=t_x}\,\ind{h\in\Omega_{\mrm{hf}}},
\end{equation}
where $Z=Z(n,\mbf{a},\mbf{b},\mbf{c},t)$ is a normalising constant, and $E$ is the set of edges of $\Z^2$ that intersect $\msf{B}_n$.
Here, `staggered' refers to the fact that the measure assigns different weights to disagreements in $h^\bullet$ and $h^\circ$ when $\mbf{a}\neq \mbf{b}$, see Section~\ref{sec:mappings_models} for details.

\subsection{Statement of main results}\label{sec:results}
We now present our main results. The proofs of each are based on a \emph{graphical representation} of the AT model originally introduced in~\cite{PfiVel97,ChaMac97} and extended in the current article, see Sections~\ref{sec:organisation}~and~\ref{sec:graph_repr}.
We emphasise that only the results based on duality are specific to the planar case. In particular, Theorems~\ref{thm:antiferro}~and~\ref{thm:sharp} and Proposition~\ref{prop:jump_mon_at} hold on $\Z^d$ for any $d\geq 2$. Define the boxes $\msf{B}_n:=\{-n,\dots,n\}^d,\,n\geq 1$.

\begin{figure}
\begin{center}
\includegraphics[page=1,height=6cm]{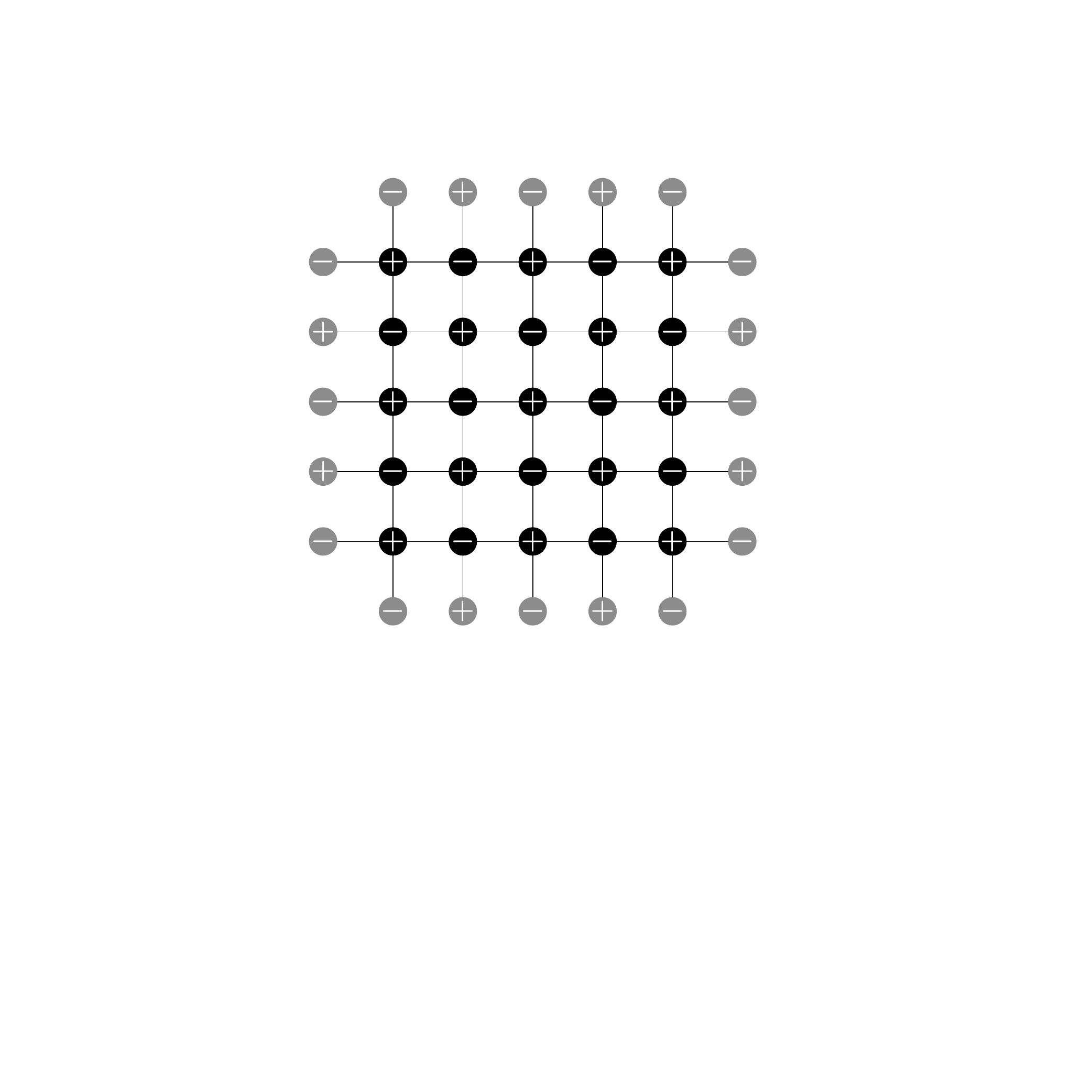}
\hspace{25pt}
\includegraphics[page=2,height=6cm]{af_at_6v.pdf}
\end{center}
\caption{\textit{Left:} The unique ground state $\eta^\pm$ of $\tau\tau'$ with respect to $\at_{\msf{B}_2,J,U}^{+,\pm}$ on the box $\msf{B}_2\subset\Z^2$ (black vertices) with external boundary (grey vertices) and edges of $\Z^2$ intersecting $\msf{B}_2$. By Theorem~\ref{thm:antiferro}, a typical configuration of $\tau\tau'$ is a small perturbation of $\eta^\pm$, whereas $\tau$ is disordered.
\textit{Right:} The unique ground state of $\hf_{\msf{B}_n,\mbf{a},\mbf{b},\mbf{c}}^{0/2,1}$ on $\msf{B}_2$ and its dual $\msf{B}_2^*\subset(\Z^2)^*$ (black hollow vertices), dual edges (dashed edges) and their endpoints outside $\msf{B}_2^*$ (grey hollow vertices). By Theorem~\ref{thm:hf_loc}, a typical height function is a small perturbation of this ground state. In relation to Theorem~\ref{thm:antiferro}, the order in the odd heights corresponds to the disorder in $\tau$, and the antiferromagnetic order in the even heights corresponds to the same property of $\tau\tau'$.}
\label{fig:af_at_6v}
\end{figure}

\paragraph{Partial antiferromagnetic phase in the AT model.}
Let $d\geq 2$. We provide the first confirmation of the presence of a phase in the AT model in which $\tau$ is \emph{disordered}, admitting exponential decay of correlations, whereas the product $\tau\tau'$ is \emph{antiferromagnetically ordered}, which is to say that correlations are bounded away from zero but have alternating signs, see Figure~\ref{fig:af_at_6v}. This was predicted in~\cite{DitBanGreKad80}, see also~\cite[Chapter 12.9]{Bax89} and~\cite{HuaDenJacSal13}.
In the corresponding part of the phase diagram, no correlation inequalities such as Griffiths' second~\cite{Gri67} or the FKG inequality~\cite{ForKasGin71} are available, and we rely on perturbative techniques.
To lighten notation and avoid case distinctions, we restrict to the isotropic case $J=J'$.

We say that a subset $\Lambda\subseteq\Z^d$ is connected if its induced subgraph is connected. Define its even and odd parts by
\[
\Lambda_{\mrm{even}}=\Big\{(x_i)_{1\leq i\leq d}\in\Lambda:\sum_{i=1}^{d}x_i\text{ even}\Big\}\quad\text{and}\quad\Lambda_{\mrm{odd}}=\Lambda\setminus\Lambda_{\mrm{even}}.
\]
Define an alternating boundary condition $\eta^{\pm}$ by setting, for $x\in\Z^d$,
\begin{equation*}
\eta^{\pm}_x:=
\begin{cases}
		+1 & \text{if }x\in\Z^d_{\mrm{even}},\\
		-1 & \text{if }x\in\Z^d_{\mrm{odd}}.
\end{cases}
\end{equation*}
We simply write $\pm$ in the superscripts of the measures.

\begin{theorem}\label{thm:antiferro}
Let $d\geq 2$. There exists $\varepsilon=\varepsilon(d)>0$ such that the following holds. For $J\in (0,\varepsilon)$ and $U\in (-\infty,-\tfrac{1}{\varepsilon})$, there exists $c=c(d,J,U)>0$ such that, for any $n\geq 1$ and any~$x\in\Z^d$,
\begin{equation}\label{eq:antiferro_1}
0\leq\langle\tau_x\rangle_{\msf{B}_n,J,U}^{+,\pm}\leq e^{-c\,(n-\norm{x})}
\qquad\text{and}\qquad
\langle\tau_x\tau'_x\rangle_{\msf{B}_n,J,U}^{+,\pm}
\begin{cases}
\geq c & \text{if }x\in\Z^d_{\mrm{even}},\\
\leq -c & \text{if }x\in\Z^d_{\mrm{odd}},
\end{cases}
\end{equation}
where $\norm{x}$ is the maximum norm of $x$.
Furthermore, there exists a random connected set $\mcal{C}\subseteq\Z^d$ satisfying $\tau_x\tau'_x=+1$ for all $x\in\mcal{C}_\mrm{even}$ and $\tau_y\tau'_y=-1$ for all $y\in\mcal{C}_\mrm{odd}$, and
\begin{equation}\label{eq:at_af_log_sized_holes}
\at_{\msf{B}_n,J,U}^{+,\pm}[\exists\,\Lambda\subseteq\Z^d\setminus\mcal{C}\text{ connected with }\mrm{diam}(\Lambda)> \log n]\leq n^{-c},
\end{equation}
where $\mrm{diam}(\Lambda)$ is the diameter of $\Lambda$ with respect to the maximum distance.
\end{theorem}

\paragraph{Ferroelectric phase in the staggered six-vertex model.}
Let $d=2$. We present an analogue of Theorem~\ref{thm:antiferro} for the corresponding six-vertex height function. It states that, within the corresponding regime, the heights of even parity are antiferromagnetically ordered, whereas the heights of odd parity exhibit ferromagnetic order, see Figure~\ref{fig:af_at_6v}. The latter is sufficient to keep the height function flat, meaning its variance is uniformly bounded in the size of the graph, which is termed \emph{localisation} of the height function.

For $n\geq 1$, denote by $\hf_{\msf{B}_n,\mbf{a},\mbf{b},\mbf{c}}^{0/2,1}$ the measure $\hf_{\msf{B}_n,\mbf{a},\mbf{b},\mbf{c}}^{t}$ when $t=2\mathbbm{1}_{\Z^2_{\mrm{odd}}}+\mathbbm{1}_{(\Z^2)^*}$, see Figure~\ref{fig:af_at_6v}.

\begin{theorem}\label{thm:hf_loc}
Let $\mbf{a},\mbf{b},\mbf{c}>0$. There exists $\varepsilon>0$ such that the following holds for $\tfrac{\mbf{a}}{\mbf{c}}<\varepsilon$ and $\tfrac{\mbf{b}}{\mbf{c}}>\tfrac{1}{\varepsilon}$. 
The variance of $h\sim\hf_{\msf{B}_n,\mbf{a},\mbf{b},\mbf{c}}^{0/2,1}$ at the origin $0\in\Z^2$ is uniformly bounded: there exists $C=C\left(\tfrac{\mbf{a}}{\mbf{c}},\tfrac{\mbf{b}}{\mbf{c}}\right)<\infty$ such that, for every $n\geq 1$,
\[
\mbf{Var}_{\msf{B}_n,\mbf{a},\mbf{b},\mbf{c}}^{0/2,1}(h_0)<C,
\]
where the variance is taken with respect to $\hf_{\msf{B}_n,\mbf{a},\mbf{b},\mbf{c}}^{0/2,1}$.
Furthermore, there exists $\alpha=\alpha\left(\tfrac{\mbf{a}}{\mbf{c}},\tfrac{\mbf{b}}{\mbf{c}}\right)>0$ and random connected sets $\mcal{C}\subseteq\Z^2$ and $\mcal{C}'\subseteq (\Z^2)^*$ such that:
\begin{enumerate}[label=(\roman*)]
\setlength\itemsep{0.1em}
\item for any $x,y\in\mcal{C}$ with $\{x,y\}\in\mbb{E}_{\Z^2},\ h_x\neq h_y$,
\item for any $x'\in\mcal{C}',\ h_{x'}=1$,
\item $\hf_{\msf{B}_n,\mbf{a},\mbf{b},\mbf{c}}^{0/2,1}[\exists\,\Lambda\subseteq\Z^2\setminus\mcal{C}\text{ or }\Lambda\subseteq(\Z^2)^*\setminus\mcal{C}'\text{ connected with }\mrm{diam}(\Lambda)> \log n]\leq n^{-\alpha}$.
\end{enumerate}
\end{theorem}

\paragraph{Subcritical sharpness along curves in the antiferromagnetic AT model.}
Let $d\geq 2$. We consider the part of the AT phase diagram where $U<0$ but where certain correlation inequalities are still valid~\cite{PfiVel97}, which corresponds to
\begin{equation}\label{eq:af_fkg_regime}\tag{AF-FKG}
\min\{J,J'\}>0>U\quad\text{and}\quad\tanh U\geq -\tanh J \tanh J'.
\end{equation}
The equation above will be referred to as both the condition and the set of $(J,J',U)\in\R^3$ that satisfy it.
One major difficulty in the analysis of this regime is the lack of general monotonicity properties in the coupling constants.
To circumvent this issue, we consider the model along smooth curves covering~\eqref{eq:af_fkg_regime}, and we show that the AT spin configuration $\tau$ undergoes a (subcritically) \emph{sharp order-disorder phase transition} along each of them.
The curves are chosen in such a way that a graphical representation of the model satisfies monotonicity properties and a suitable Russo-type inequality along them, allowing to run the general argument of~\cite{DumRaoTas19}.

Given a parametrised curve $\gamma:(0,1)\to\R^3$ and $\beta\in (0,1)$, we write $\langle\,\cdot\,\rangle_{\Lambda,\gamma(\beta)}^{\eta,\eta'}$ for the expectation operator with respect to the measure in~\eqref{eq:def_at} with coupling constants given by $\gamma(\beta)$.

\begin{theorem}\label{thm:sharp}
There exists a family of disjoint smooth curves $\gamma_{\kappa,\kappa'}:(0,1)\to\R^3,$ ${0<\kappa<\kappa'\leq 1,}$ that covers~\eqref{eq:af_fkg_regime} and that satisfies the following.
For $d\geq 2$ and any $\kappa,\kappa'\in(0,1]$ with $\kappa<\kappa'$, there exists $\beta_\mrm{c}=\beta_\mrm{c}(d,\kappa,\kappa')\in (0,1)$ such that
\begin{itemize}
\setlength\itemsep{0.1em}
\item for $\beta<\beta_\mrm{c}$, there exists $c_\beta=c_\beta(d,\kappa,\kappa')>0$ such that $0\leq\langle\tau_0\rangle_{\msf{B}_n,\gamma_{\kappa,\kappa'}(\beta)}^{+,\mrm{f}}\leq e^{-c_\beta n}$,
\item there exists $c=c(d,\kappa,\kappa')>0$ such that, for $\beta>\beta_\mrm{c}$, $\langle\tau_0\rangle_{\msf{B}_n,\gamma_{\kappa,\kappa'}(\beta)}^{+,\mrm{f}}\geq c(\beta-\beta_\mrm{c})$.
\end{itemize}
\end{theorem}

The above theorem provides a foundation for validating an analogous statement in the isotropic model when $d=2$.
In fact, once monotonicity properties that we conjecture are confirmed, a sharp phase transition at the self-dual curve $\sinh 2J=e^{-2U}$~\cite{MitSte71} in the isotropic model is a consequence of Theorem~\ref{thm:sharp}.
See Section~\ref{sec:sharp_af_iso}, where we also show indications of monotonicity in the isotropic regime, which imply, for example, the following. 

\begin{figure}
\begin{center}
\includegraphics[height=6.5cm]{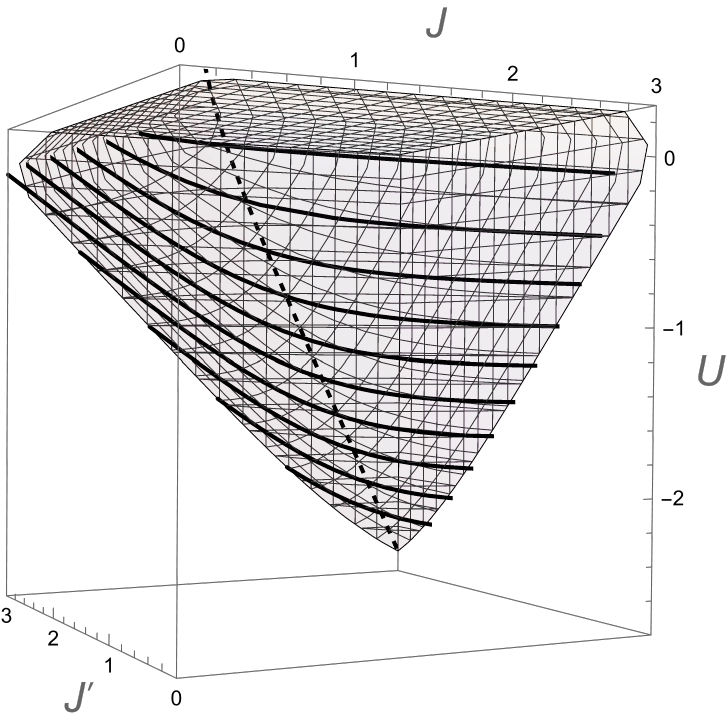}
\quad
\includegraphics[height=6.5cm]{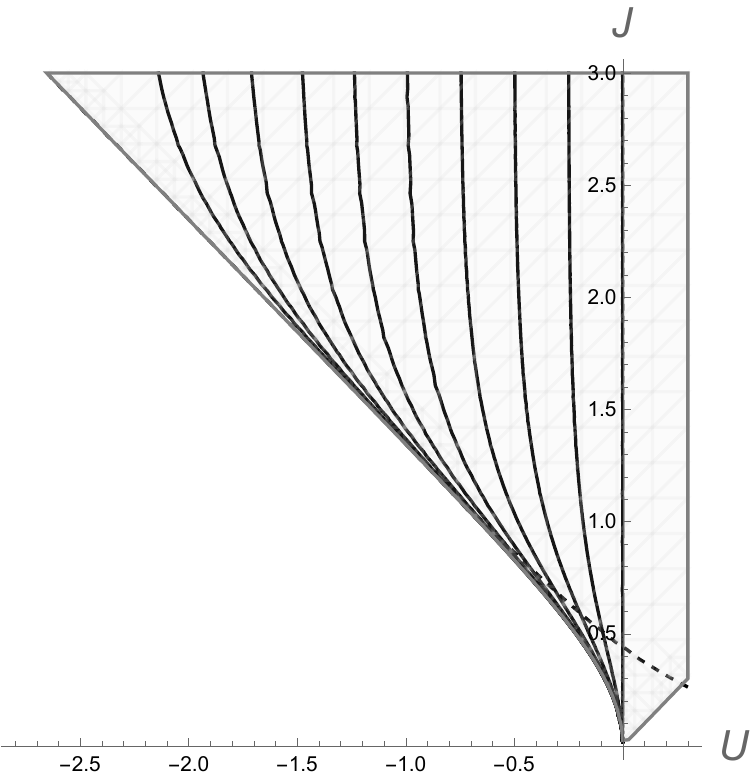}
\end{center}
\caption{\textit{Left:} Part of the regime~\eqref{eq:af_fkg_regime} (light grey area) and curves $\gamma_{\kappa,\kappa'}$ (solid black) that intersect the self-dual isotropic curve $\sinh 2J=e^{-2U}$ (dashed black), for ten values of $\kappa,\,\kappa'$. \textit{Right:} Part of the regime~\eqref{eq:af_fkg_regime} (light grey area) in the isotropic phase diagram and curves $\widehat{\gamma}_\kappa$ (solid black) that intersect the self-dual curve (dashed black), for ten values of $\kappa$.}
\label{fig:curves}
\end{figure}

\begin{proposition}\label{prop:jump_mon_at}
There exists a family of disjoint smooth curves $\widehat{\gamma}_\kappa:(0,1)\to\R^3,\,\kappa\in (0,1),$ covering the isotropic part of~\eqref{eq:af_fkg_regime}, and smooth, strictly increasing, bijective functions ${f_\kappa:(0,1)\to (0,\infty)},\,\kappa\in (0,1),$ such that the following holds for $d\geq 2$. 
For any $\kappa\in (0,1)$, any $\beta_1,\beta_2\in (0,1)$ with $f_\kappa(\beta_1) \leq \kappa\,f_\kappa(\beta_2)$, and any $\Lambda\subset\Z^d$, 
\begin{equation}
\langle\tau_0\rangle_{\Lambda,\widehat{\gamma}_{\kappa}(\beta_1)}^{+,\mrm{f}}\leq \langle\tau_0\rangle_{\Lambda,\widehat{\gamma}_{\kappa}(\beta_2)}^{+,\mrm{f}}.
\end{equation}
\end{proposition}
We also prove that, along the same curves $\widehat{\gamma}_\kappa$, a sum of infinite-volume edge densities of the graphical representation is non-decreasing, see Proposition~\ref{prop:mon_gat_edge_densities}. See Figure~\ref{fig:curves} for an illustration of the regime~\eqref{eq:af_fkg_regime} and the curves $\gamma_{\kappa,\kappa'}$ and $\widehat{\gamma}_\kappa$.

\subsection{Overview and organisation of the article}\label{sec:organisation}
The article is organised as follows.
\begin{itemize}
\setlength\itemsep{0.1em}
\item Section~\ref{sec:setting} is devoted to the introduction and clarification of some notation and conventions.

\item In Section~\ref{sec:graph_repr}, we revisit a \emph{graphical representation} (GAT) which, in certain parts of the phase diagram, coincides with those introduced in~\cite{PfiVel97,ChaMac97}. 
This is a \emph{percolation} model obtained by means of an FK-Ising-type expansion for one of the AT spin configurations.
In contrast to~\cite{PfiVel97}, where a pair of percolation configurations is considered, we restrict ourselves to one of them to enable greater generality.
An immediate consequence is the validation of \emph{Griffiths' first inequality}~\cite{Gri67} for the AT spin configuration $\tau$ when $J\geq\min\{\abs{J'},\abs{U}\}$, see Corollary~\ref{cor:gks}.
Subsequently, we discuss the property of \emph{positive association} of the graphical representation and its relation to~\cite{PfiVel97} by constructing a coupling of a pair of them.

\item In Section~\ref{sec:af_phase}, we utilise the graphical representation to prove Theorem~\ref{thm:antiferro}. 

\item In Section~\ref{sec:mappings_models}, we make use of the graphical representation again to formulate the connection between the AT and the eight-vertex models by means of coupling. We then discuss the special case of the six-vertex model and the relation of Theorems~\ref{thm:antiferro} and~\ref{thm:hf_loc}, and we present the proof of Theorem~\ref{thm:hf_loc}.

\item In Section~\ref{sec:sharpness}, we discuss the concept of subcritically sharp order-disorder phase transitions in the AT model. We then use the graphical representation once more to construct the curves from Theorem~\ref{thm:sharp} satisfying the desired properties, and we prove the theorem. Following a brief discussion of the isotropic case, we propose a conjecture concerning monotonicity. The confirmation of this conjecture would enable us to deduce the presence of a sharp order-disorder phase transition in the isotropic model based on our findings for the anisotropic case. Subsequently, we substantiate this conjecture by presenting indications of monotonicity within the isotropic model, see Propositions~\ref{prop:mon_gat_edge_densities}~and~\ref{prop:jump_mon_gat}.

\item Appendices~\ref{sec:proof_fkg_lattice}--\ref{sec:jump_mon_gat} contain the proof of the FKG-lattice condition for the graphical representation (Lemma~\ref{lem:gat_fkg}), the proof of maximality of certain boundary conditions (Lemma~\ref{lem:max_bc}), the construction of the curves for Theorem~\ref{thm:sharp} (Lemma~\ref{lem:curves}), and the proof of a `jump monotonicity' statement in the antiferromagnetic isotropic regime (Proposition~\ref{prop:jump_mon_gat}). 
\end{itemize}

\paragraph{Acknowledgements.} I am particularly grateful to Alexander Glazman for his encouragement in the further investigation of some of the issues that are the subject of this article, and for useful comments on a first draft. I also thank him, as well as Tom{\'a}s Alcalde, Yacine Aoun, Marcin Lis, S{\'e}bastien Ott and Kieran Ryan for stimulating suggestions and discussions on the Ashkin--Teller and six-vertex models.
This research was funded in whole or in part by the Austrian Science Fund (FWF) [10.55776/P34713].
The author thanks the anonymous referees and associate editor for their constructive comments.

\section{Notation and conventions}\label{sec:setting}

\paragraph{Graph notation.}
Let $\mbb{G}=(\mbb{V}_\mbb{G},\mbb{E}_\mbb{G})$ be a graph.
For an edge $\{x,y\}\in\mbb{E}_{\mbb{G}}$, we simply write $xy=\{x,y\}$.
Given $V\subset \mbb{V}_\mbb{G}$, define
\begin{align*}
\overline{V}&:=V\cup\{y\in \mbb{V}_\mbb{G}:\exists x\in V,\,xy\in\mbb{E}_\mbb{G}\},\\
\overline{\mbb{E}}_V&:=\{e\in\mbb{E}_\mbb{G} :e\cap V\neq\varnothing\}.
\end{align*}

A \emph{path} in $\mbb{G}$ is a sequence $\Gamma=(x_0,\dots,x_{\ell})$ of pairwise distinct vertices $x_i\in\mbb{V}_\mbb{G}$ such that $x_ix_{i+1}\in\mbb{E}_\mbb{G}$ for $i=0,\dots,\ell-1$.
We call $x_0x_1,\dots,x_{\ell-1}x_\ell$ the edges of $\Gamma$, and $\ell=:\abs{\Gamma}$ the length of $\Gamma$. When $x_0=x_\ell$, we call $\Gamma$ a \emph{circuit}.

The integer lattice has vertex-set $\Z^d=\{(x_i)_{1\leq i\leq d}\in\R^d:x_i\in\Z\text{ for all }i\}$ and edges between nearest neighbours, that is, between vertices of Euclidean distance 1. In a slight abuse of notation, we will also write $\Z^d$ for the graph itself.
For $n\geq 0$, define $\mathsf{B}_n:=\{-n,\dots,n\}^d\subset\Z^d$ and $\msf{B}_{-1}=\varnothing$, and, for $x\in\Z^d$, let $\msf{B}_n(x)$ be the translate of $\msf{B}_n$ by $x$. 

\paragraph{Spin and percolation configurations.}
Let $\mbb{G}=(\mbb{V}_\mbb{G},\mbb{E}_\mbb{G})$ be a graph. Given $V\subseteq\mbb{V}_\mbb{G}$, we call $\spin=(\spin_x)_{x\in V}\in\{-1,1\}^{V}$ a \emph{spin configuration} on $V$, and we write $\Sigma_{V}$ for the set of all such configurations.
When $V=\mbb{V}_\mbb{G}$, we simply write $\Sigma_\mbb{G}$ or $\Sigma$.
Moreover, given $\eta\in\Sigma_{\mbb{G}}$, define 
\begin{align*}
\Sigma_V^\eta&:=\{\spin\in\Sigma_{\mbb{G}}:\spin_x=\eta_x \text{ for all }x\in \mbb{V}_\mbb{G}\setminus V\},\\
\Sigma_V^{\mrm{f}}&:=\{\spin\in\Sigma_{\mbb{G}}:\spin_x=1 \text{ for all }x\in \mbb{V}_\mbb{G}\setminus \overline{V}\}.
\end{align*}

\begin{remark}\label{rem:free_bc}
Observe that this definition of free boundary conditions is slightly different from the standard one in the Ising model, and in fact coincides with the latter on the graph~\((\overline{V},\overline{\mbb{E}}_V)\).
The reason for this choice is that the interactions for free boundary conditions and for those given by a configuration~\(\eta\) will then occur along the same edges, thus simplifying the presentation and avoiding case distinctions.
\end{remark}

The restriction $(\spin_x)_{x\in V}$ of $\spin$ to $V$ is denoted by $\spin_V$.
For a pair $\spin,\spin'\in\Sigma_V$, we write $\spin\spin':=(\spin_x \spin'_x)_{x\in V}$ for their coordinatewise product.
Given $E\subseteq\mbb{E}_{\mbb{G}}$, define $E_\spin$ as the set of \emph{disagreement edges} of $\spin$ in $E$, that is,
\begin{equation}\label{eq:def_disagr_edges}
E_\spin:=\{xy\in E: \spin_x\neq \spin_y\}.
\end{equation}

We call $\omega=(\omega_e)_{e\in E}\in\{0,1\}^{E}$ a \emph{percolation configuration} on $E$, and we write $\Omega_E$ for the set of all such configurations. An edge $e\in E$ is said to be \emph{open} (respectively, \emph{closed}) in $\omega$ if $\omega_e=1$ (respectively, $\omega_e=0$), and a path or a set of edges is open (respectively, closed) if all its edges are open (respectively, closed).
We identify $\omega$ with the set of open edges $\{e\in E:\omega_e=1\}$ as well as with the subgraph $(\mbb{V}_\mbb{G},\omega)$. The connected components of $\omega$ are called \emph{clusters}. When $E=\mbb{E}_{\mbb{G}}$, we simply write $\Omega_\mbb{G}$ or $\Omega$. Moreover, given $\omega\in\Omega_{\mbb{G}}$, we write $\Omega_E^\omega$ for the set of percolation configurations on $\mbb{E}_\mbb{G}$ that coincide with $\omega$ on $\mbb{E}_\mbb{G}\setminus E$. The restriction $(\omega_e)_{e\in E}$ of $\omega$ to $E$ is denoted by $\omega_E$. Given $V\subset\mbb{V}_\mbb{G}$, we denote by $k_V(\omega)$ the number of clusters of $\omega$ that intersect $\overline{V}$. For $V,W\subset\mbb{V}_\mbb{G}$, we write $V\xleftrightarrow{\omega} W$ (or simply $V\leftrightarrow W$ when no confusion is possible) if there exists an open path in $\omega$ that connects a vertex in $V$ to a vertex in $W$, and we write $V\xleftrightarrow{\omega}\infty$ (or simply $V\leftrightarrow\infty$) if some vertex in $V$ belongs to an infinite cluster in $\omega$. In this context, singleton sets are denoted without their braces.

The sets of configurations $\Sigma_V$ and $\Omega_E$ as well as products of them are equipped with the natural product $\sigma$-algebra.

\paragraph{Conventions.}
We use bold letters $\mbf{P,\,E,\,Var,\,Cov}$ to denote probability measures, expectations, variances and covariances, respectively. 
To lighten notation, we will use the symbols $\tau,\spin,\sigma,\omega$ for both deterministic configurations and random variables. More precisely, when we consider a random variable $X$, we may express its law by denoting a deterministic element in its target space by the same symbol $X$. 
For a probability measure $\mu$ on a measurable space $(M,\mcal{M})$ and $m\in M$, we write $\mu[m]$ instead of $\mu[\{m\}]$. Given in addition a function $f:M\to [0,\infty)$, we write $\mu[m]\propto f(m),\, m\in M,$ if $\mu$ is proportional to $f$, that is, if there exists a constant $C>0$ such that $\mu[m]=C f(m)$ for all $m\in M$.

\section{Graphical representation of the AT model}\label{sec:graph_repr}
\subsection{The representation}\label{sec:graph_repr_def}
As in the classical Edwards--Sokal coupling~\cite{EdwSok88}, we are going to define a percolation model whose connection probabilities describe correlations of one of the AT spin configurations. It is obtained by performing an FK-Ising-type expansion for one of the AT spin configurations while keeping the other fixed.
For certain ranges of the coupling constants, this representation was already considered in~\cite{PfiVel97,ChaMac97,Lis22,RaySpi22,GlaPel23,AouDobGla24}, see Remark~\ref{rem:gat_ref}. As the spin configuration may be either $\tau,\,\tau'$ or the product $\tau\tau'$ and we formulate the coupling only once, we shall apply a change of variables later on, which reorders the coupling constants $J,J',U$. For instance, if $(\tau,\tau')$ is distributed according to $\at_{\Lambda,J,J',U}^{+,+}$, then $(\tau,\tau\tau')$ is distributed according to $\at_{\Lambda,J,U,J'}^{+,+}$. For this purpose, we will use different symbols $\spin,\spin'$ for spin configurations and different letters $K,K',K''\in\R$ for the coupling constants. Define
\begin{equation}\label{eq:sample_edge_densities}
p_1=1-e^{-2(K-K'')}\qquad\text{and}\qquad p_2=1-e^{-2(K+K'')},
\end{equation}
and, in order for them to be densities, we will make the assumption that
\[
K\geq\abs{K''}.
\]
Let $\Lambda\subset\Z^d$ be finite, and set $E=\overline{\mbb{E}}_\Lambda$.
Fix boundary conditions $\eta\in\{+,\mrm{f}\}$ and $\eta'\in\Sigma_{\Z^d}\cup\{\mrm{f}\}$ arbitrary, and set $\#=\ind{\eta=+}$.
Take $(\spin,\spin')$ distributed according to $\at_{\Lambda,K,K',K''}^{\eta,\eta'}$, and sample a percolation configuration $\omega\in\Omega_{\Z^d}$ as follows. 
Recall the definition~\eqref{eq:def_disagr_edges} of $E_\spin$.
Define $\omega_e$ independently for each edge $e\in\mbb{E}_{\Z^d}$: if $e\notin E$, set $\omega_e=\#$. Otherwise, set $\omega_e=0$ for $e\in E_\spin$, $\omega_e=1$ with probability $p_1$ for $e\in E_{\spin'}\setminus E_\spin$, and $\omega_e=1$ with probability $p_2$ for $e\in E\setminus(E_{\spin}\cup E_{\spin'})$. Formally, let $(U_e)_{e\in E}$ be i.i.d. uniform on $[0,1]$, independent of $(\spin,\spin')$, and define 
\begin{equation}\label{eq:sample_rule_gat}
\omega_e=
	\begin{cases}
		\# & \text{if }e\notin E,\\
		0 & \text{if }e\in E_\spin,\\
		\bbm{1}_{[0,p_1]}(U_e) & \text{if }e\in E_{\spin'}\setminus E_\spin,\\
		\bbm{1}_{[0,p_2]}(U_e) & \text{if }e\in E\setminus(E_{\spin}\cup E_{\spin'}).
	\end{cases}
\end{equation}
We denote the law of $\omega$ on $\Omega_{\Z^d}$ by $\gat_{\Lambda,K,K',K''}^{\#,\eta'}$. Recall that, by convention, we use the same symbols for both random variables and their deterministic realisations. 

\begin{proposition}\label{prop:at_gat_laws}
\noindent(i) The joint law of $(\spin,\spin',\omega)$ on $(\Sigma_{\Z^d})^2\times\Omega_{\Z^d}$ is given by 
\begin{equation}\label{eq:at_gat_joint_law}
\begin{aligned}
\Prob{\spin,\spin',\omega}\propto e^{2(K''-K')\abs{E_{\spin'}}}\,\big(e^{2(K-K'')}-1\big)^{\abs{\omega_E\cap E_{\spin'}}}\,\big(e^{2(K+K'')}-1\big)^{\abs{\omega_E\setminus E_{\spin'}}}\\
\cdot\,\ind{(\spin,\spin')\in\Sigma_\Lambda^\eta\times\Sigma_\Lambda^{\eta'}}\,\ind{\omega\in\Omega_E^\#}\,\ind{\omega\cap E_\spin=\varnothing}.
\end{aligned}
\end{equation}
In particular, conditionally on $\omega$, the law of $\spin_\Lambda$ is obtained as follows. If $\eta=\mrm{f}$, assign $\pm1$ uniformly and independently to clusters of $\omega$ that intersect $\Lambda$. If $\eta=+$, assign $+1$ to clusters of $\omega$ that intersect $\Z^d\setminus\Lambda$ and $\pm1$ uniformly and independently to the other ones.\vspace{5pt}

\noindent(ii) When $K+K''>0$, the law of $\omega$ on $\Omega_{\Z^d}$ is given by 
\begin{equation}\label{eq:gat_law_1}
\gat_{\Lambda,K,K',K''}^{\#,\eta'}[\omega]\propto \mrm{w}_1^{\abs{\omega_E}}\,2^{k_{\Lambda}(\omega)}\,\ind{\omega\in\Omega_E^\#}\sum_{\spin'\in\Sigma_\Lambda^{\eta'}}\mrm{w}_2^{\abs{E_{\spin'}\setminus\omega}}\,\mrm{w}_3^{\abs{E_{\spin'}\cap\omega}},
\end{equation}
where the weights are defined as
\begin{equation}\label{eq:gat_weights_1}
\mrm{w}_1:=e^{2(K+K'')}-1,\quad\mrm{w}_2:=e^{2(K''-K')}\quad\text{and}\quad \mrm{w}_3:=\frac{e^{-2K'}\big(e^{2K}-e^{2K''}\big)}{e^{2(K+K'')}-1}.
\end{equation}

\noindent(iii) When $K+K''=0$, the law of $\omega$ on $\Omega_{\Z^d}$ is given by 
\begin{equation*}
\gat_{\Lambda,K,K',K''}^{\#,\eta'}[\omega]\propto \big(e^{4K}-1\big)^{\abs{\omega_E}}\,2^{k_{\Lambda}(\omega)}\,\ind{\omega\in\Omega_E^\#}\sum_{\substack{\spin'\in\Sigma_\Lambda^{\eta'}:\\ E_{\spin'}\supseteq\omega_E}}e^{-2(K+K')\abs{E_{\spin'}}}.
\end{equation*}
\end{proposition}

As a consequence of part (i) of the above proposition (see, e.g.,~\cite[Theorem 1.16]{Gri06}), correlations of $\spin$ are described by connection probabilities of $\omega$.
\begin{corollary}\label{cor:gat_corr_conn}
For any $x,y\in\Lambda$,
\[
\langle \spin_x\spin_y \rangle_{\Lambda,K,K',K''}^{\eta,\eta'}=\gat_{\Lambda,K,K',K''}^{\#,\eta'}[x\leftrightarrow y]\quad\text{and}\quad
\langle \spin_x \rangle_{\Lambda,K,K',K''}^{+,\eta'}=\gat_{\Lambda,K,K',K''}^{1,\eta'}[x\leftrightarrow\Z^d\setminus\Lambda].
\]
\end{corollary}

The above proposition also provides an example of non-uniqueness of graphical representations. For instance, consider the AT measure in~\eqref{eq:def_at} with $J=J'>U>0$. Then one can choose both $(\spin,\spin')=(\tau,\tau')$ and $(\spin,\spin')=(\tau,\tau\tau')$ to obtain two different percolation models for which the connection probabilities describe correlations of the AT spin configuration $\tau$. Only the former enjoys the property of positive association, see Section~\ref{sec:pos_ass}.

\begin{proof}[Proof of Proposition~\ref{prop:at_gat_laws}]
\noindent(i) Write $\spin_x\spin_y=1-2\ind{\spin_x\neq \spin_y}$ and analogously for $\spin'_x\spin'_y$ to obtain that 
\begin{equation}\label{eq:at_law_2}
\at_{\Lambda,K,K',K''}^{\eta,\eta'}[\spin,\spin']\propto e^{-2(K+K'')\abs{E_\spin}}\,e^{-2(K'+K'')\abs{E_{\spin'}}}\,e^{4K''\abs{E_\spin\cap E_{\spin'}}}\,\ind{(\spin,\spin')\in\Sigma_\Lambda^\eta\times\Sigma_\Lambda^{\eta'}}.
\end{equation}
Moreover, by~\eqref{eq:sample_rule_gat}, we have
\begin{alignat*}{2}
\Prob{\spin,\spin',\omega}&=\at_{\Lambda,K,K',K''}^{\eta,\eta'}[\spin,\spin']&&\cdot\,p_1^{\,\abs{\omega_E\cap E_{\spin'}}}\,(1-p_1)^{\abs{E_{\spin'}\setminus(\omega_E\cup E_\spin)}}\\
& &&\cdot\,p_2^{\,\abs{\omega_E\setminus E_{\spin'}}}\,(1-p_2)^{\abs{E\setminus (\omega_E\cup E_\spin\cup E_{\spin'})}}\,\ind{\omega\in\Omega_E^\#}\,\ind{\omega\cap E_\spin=\varnothing}\\
&=\at_{\Lambda,K,K',K''}^{\eta,\eta'}[\spin,\spin']&&\cdot\,p_1^{\,\abs{\omega_E\cap E_{\spin'}}}\,(1-p_1)^{\abs{E_{\spin'}\setminus E_\spin}-\abs{\omega_E\cap E_{\spin'}}}\\
& &&\cdot\,p_2^{\,\abs{\omega_E\setminus E_{\spin'}}}\,(1-p_2)^{\abs{E}-\abs{E_{\spin}}-\abs{E_{\spin'}\setminus E_\spin}-\abs{\omega_E\setminus E_{\spin'}}}\\
& &&\cdot\,\ind{\omega\in\Omega_E^\#}\,\ind{\omega\cap E_\spin=\varnothing},
\end{alignat*}
where we have repeatedly used the fact that $\omega_E$ and $E_\spin$ are disjoint, which is guaranteed by the corresponding indicator. 
Plugging in~\eqref{eq:sample_edge_densities} and~\eqref{eq:at_law_2} and collecting like terms, we arrive at~\eqref{eq:at_gat_joint_law}.

\noindent(ii) Writing $\abs{\omega_E\setminus E_{\spin'}}=\abs{\omega_E}-\abs{\omega_E\cap E_{\spin'}}$ and collecting like terms, equation~\eqref{eq:at_gat_joint_law} turns into
\begin{equation*}
\Prob{\spin,\spin',\omega}\propto \mrm{w}_1^{\abs{\omega_E}}\,\mrm{w}_2^{\abs{E_{\spin'}\setminus\omega}}\,\mrm{w}_3^{\abs{E_{\spin'}\cap\omega}}\,\ind{(\spin,\spin')\in\Sigma_\Lambda^\eta\times\Sigma_\Lambda^{\eta'}}\,\ind{\omega\in\Omega_E^\#}\,\ind{\omega\cap E_\spin=\varnothing}.
\end{equation*}
Summing over $(\spin,\spin')$ while noting that there exist $2^{k_\Lambda(\omega)+\mrm{const}(\Lambda,\#)}$ spin configurations $\spin\in\Sigma_\Lambda^\eta$ with $E_\spin\cap\omega=\varnothing$, we obtain~\eqref{eq:gat_law_1}.

\noindent(iii) When $K+K''=0$, equation~\eqref{eq:at_gat_joint_law} turns into 
\[
e^{-2(K+K')\abs{E_{\spin'}}}\,\big(e^{4K}-1\big)^{\abs{\omega_E}}\,\ind{\omega_E\subseteq E_{\spin'}}\,\ind{(\spin,\spin')\in\Sigma_\Lambda^\eta\times\Sigma_\Lambda^{\eta'}}\,\ind{\omega\in\Omega_E^\#}\,\ind{\omega\cap E_\spin=\varnothing},
\]
and summing over $(\spin,\spin')$ as in (ii) finishes the proof.
\end{proof}

\begin{remark}\label{rem:gat_ref}
Let $\omega$ be distributed according to $\gat_{\Lambda,K,K',K''}^{\#,\eta'}$.
The following equalities of laws hold up to boundary conditions and for $\Lambda$ chosen accordingly.
\begin{enumerate}[label=(\roman*)]
\setlength\itemsep{0.1em}
\item When $\min\{K,K'\}\geq \max\{K'',0\}$ and $\tanh K''\geq -\tanh K \tanh K'$ (which implies $K\geq \abs{K''}$), the law of $\omega$ coincides with the law of $\underline{n}_\sigma$ under $\nu_\Lambda^{+}(\,\cdot\,\vert 2,2)$ with $J_\sigma=K,\,J_\tau=K'$ and $J_{\sigma\tau}=K''$ in~\cite[Proposition 3.1]{PfiVel97}.
\item When $K>K'=K''$ and $e^{-2K}=\sinh 2K'$, the law of $\omega_E$ coincides with~\cite[Equation (7.2)]{GlaPel23} with $c=\coth 2K'$. Moreover, when $K=K''$ and $\sinh 2K=e^{-2K'}$, it coincides with the law of $\xi^*$ in~\cite[Lemmata 7.1 and 8.1]{GlaPel23} with $c=\coth 2K$.
\item When $K>K'=K''>0$, the law of $\omega_E$ coincides with~\cite[Equation (3.6)]{Lis22} for $a=\tanh 2K'$ and $b=(e^{2K} \cosh 2K')^{-1}$, where it appears in the context of the staggered six-vertex model, see Section~\ref{sec:mappings_models}.
\item When $K>K'=K''$ and $e^{-2K}=\sinh 2K'$, the law of $\omega$ coincides with the marginal of~\cite[Equation (5.10)]{RaySpi22} on $\eta^0$ with $\alpha=\coth (2K')-2$ and $q=2$.
\item When $K=K''\geq K'$ and $\sinh 2K=e^{-2K'}$, the law of $\omega_E$ coincides with the law of $\xi^\bullet$ for $a=b=1$ and $c=\coth 2K\leq 2$ in~\cite[Section 4]{GlaLam23}, where it is decomposed into a sum to recover positive association, which $\omega$ itself does not satisfy in that regime, see Section~\ref{sec:pos_ass}.
\end{enumerate}
\end{remark}

The graphical representation satisfies a strong \emph{finite-energy} property which will be crucial for the perturbative analyses in Sections~\ref{sec:af_phase} and~\ref{sec:8v_ferroel_phase}.

\begin{lemma}\label{lem:finite_energy}
Let $K+K''>0$. For any $e\in E$ and any $\xi\in\Omega_E^{\#}$ with $\gat_{\Lambda,K,K',K''}^{\#,\eta'}[\xi]>0$,
\[
\Big(1+\mrm{w}_1^{-1}2\,\tfrac{\max\{1,\mrm{w}_2\}}{\min\{1,\mrm{w}_3\}}\Big)^{-1}
\leq
\gat_{\Lambda,K,K',K''}^{\#,\eta'}\big[\omega_e=1\,\big|\,\omega_{\,\mbb{E}_{\Z^d}\setminus\{e\}}=\xi_{\,\mbb{E}_{\Z^d}\setminus\{e\}}\big]
\leq 
\Big(1+\mrm{w}_1^{-1}\,\tfrac{\min\{1,\mrm{w}_2\}}{\max\{1,\mrm{w}_3\}}\Big)^{-1},
\]
where the weights $\mrm{w}_i$ are given by~\eqref{eq:gat_weights_1}, and where the quantity on the left side is defined to be $0$ when $\mrm{w}_3=0$, that is, when $K=K''$. 
\end{lemma}
It should be mentioned that, for the case $\mrm{w}_3=0$, one can also obtain a non-trivial lower bound, which is omitted as it is not used in this article.
\begin{proof}
Let $F=E\setminus\{e\}$. Equation~\eqref{eq:gat_law_1} implies
\begin{align*}
\gat&_{\Lambda,K,K',K''}^{\#,\eta'}\big[\omega_e=1\,\big|\,\omega_{\,\mbb{E}_{\Z^d}\setminus\{e\}}=\xi_{\,\mbb{E}_{\Z^d}\setminus\{e\}}\big]\\
&=\left(1+\mrm{w}_1^{-1}\,2^{k_{\Lambda}(\xi\setminus\{e\})-k_{\Lambda}(\xi\cup\{e\})}\frac{\sum_{\spin'}\mrm{w}_2^{\abs{E_{\spin'}\setminus (\xi\setminus\{e\})}}\mrm{w}_3^{\abs{E_{\spin'}\cap (\xi\setminus\{e\})}}}{\sum_{\spin'}\mrm{w}_2^{\abs{E_{\spin'}\setminus (\xi\cup\{e\})}}\mrm{w}_3^{\abs{E_{\spin'}\cap (\xi\cup\{e\})}}}\right)^{-1}\\
&=\left(1+\mrm{w}_1^{-1}\,2^{k_{\Lambda}(\xi\setminus\{e\})-k_{\Lambda}(\xi\cup\{e\})}\frac{\sum_{\spin'}\mrm{w}_2^{\abs{F_{\spin'}\setminus \xi}}\mrm{w}_3^{\abs{F_{\spin'}\cap\xi}}\mrm{w}_2^{\ind{e\in E_{\spin'}}}}{\sum_{\spin'}\mrm{w}_2^{\abs{F_{\spin'}\setminus \xi}}\mrm{w}_3^{\abs{F_{\spin'}\cap\xi}}\mrm{w}_3^{\ind{e\in E_{\spin'}}}}\right)^{-1},
\end{align*}
where the summations are over all $\spin'\in\Sigma_\Lambda^{\eta'}$.
The exponent $k_{\Lambda}(\xi\setminus\{e\})-k_{\Lambda}(\xi\cup\{e\})$ is either $0$ or $1$, and the claim of the lemma follows. 
\end{proof}

\subsection{Griffiths' first inequality}\label{sec:griffiths}
For this subsection, we return to the terminology of~\eqref{eq:def_at}. As a direct consequence of Proposition~\ref{prop:at_gat_laws}, we obtain Griffiths' first inequality~\cite{Gri67} for $\tau$ when $J\geq\min\{\abs{J'},\abs{U}\}$, which was previously established when $J,J',U\geq 0$ \cite{KelShe68} and more generally when $\min\{J,J'\}\geq \max\{U,0\}$ and $\tanh U\geq -\tanh J \tanh J'$ or the same holds with $J'$ and $U$ interchanged~\cite{PfiVel97}. 

For $\tau\in\Sigma_{\Z^d}$ and $A\subset\Z^d$ finite, define 
\begin{equation*}
\tau^A:=\prod_{x\in A}\tau_x.
\end{equation*}

\begin{corollary}\label{cor:gks}
Let $J,J',U\in\R$ with $J\geq \min\{\abs{J'},\abs{U}\}$. Then, for $A\subseteq\Lambda\subset\Z^d$ finite,
\begin{equation}\label{eq:gks}
\langle \tau^A \rangle^{+,+}_{\Lambda,J,J',U}\geq 0.
\end{equation} 
\end{corollary}
\begin{proof}
It is classical (see, e.g.,~\cite[Proposition 3.1]{PfiVel97} and its proof) that Proposition~\ref{prop:at_gat_laws}(i) implies that, for $K,K',K''$ with $K\geq\abs{K''}$,
\[
\langle \spin^A \rangle^{+,+}_{\Lambda,K,K',K''}=\gat_{\Lambda,K,K',K''}^{1,+}[\kappa_A]\geq 0,
\]
where $\kappa_A$ is the event that every finite cluster of $\omega$ contains an even number of vertices in $A$. Applying this for $\spin=\tau,\,\spin'=\tau'$ and $(K,K',K'')=(J,J',U)$ gives~\eqref{eq:gks} for $J\geq\abs{U}$, while the choice $\spin=\tau,\,\spin'=\tau\tau'$ and $(K,K',K'')=(J,U,J')$ gives~\eqref{eq:gks} for $J\geq\abs{J'}$.
\end{proof}

\begin{remark}
\begin{enumerate}[label=(\roman*)]
\item For the isotropic AT model when $J=J'$, the condition in Corollary~\ref{cor:gks} reduces to $J\geq 0$, and is therefore optimal.
\item Another way to prove Corollary~\ref{cor:gks} is to consider $\tau$ as an Ising model with random coupling constants and then apply~\cite{KelShe68} to the conditional measures. 
\end{enumerate}
\end{remark}

\subsection{Positive association}\label{sec:pos_ass}
In some part of the phase diagram, the graphical representation is \emph{positively associated}. We point out that we require this property only in Section~\ref{sec:sharpness}.
Fix a graph $\mbb{G}=(\mbb{V}_\mbb{G},\mbb{E}_\mbb{G})$ and $E\subseteq\mbb{E}_\mbb{G}$ finite. We equip the set of percolation configurations $\Omega_E$ with a partial order $\leq$ which is defined coordinatewise: for $\omega,\omega'\in\Omega_E,\ \omega\leq\omega'$ if and only if $\omega_e\leq\omega'_e$ for all $e\in E$. 
A random variable $X:\Omega_E\to\R$ is called \emph{increasing} if $X(\omega)\leq X(\omega')$ whenever $\omega\leq\omega'$.
An event $A\subseteq\Omega_E$ is called increasing if its indicator $\bbm{1}_A$ is increasing.
Given two measures $\mu,\,\nu$ on $\Omega_E$, we say that $\mu$  is (stochastically) \emph{dominated} by $\nu$ (or $\nu$ dominates $\mu$), and we write $\mu\leq_{\mrm{st}}\nu$ (or $\nu\geq_{\mrm{st}}\mu$), if $\mbf{E}_\mu[X]\leq \mbf{E}_\nu[X]$ for all increasing $X:\Omega_E\to[0,\infty)$.
Furthermore, we say that a measure $\mu$ on $\Omega_E$ is \emph{positively associated} if it satisfies the \emph{FKG inequality}:
\begin{equation}\label{eq:fkg}\tag{FKG}
\mbf{E}_\mu [XY]\geq \mbf{E}_\mu [X]\mbf{E}_\mu [Y]\qquad\text{for all increasing }X,Y:\Omega_E\to[0,\infty).
\end{equation}
Moreover, $\mu$ is said to satisfy the \emph{FKG lattice condition} if
\begin{equation}\label{eq:fkg_lattice}\tag{FKG-L}
\mu [\omega\vee\omega']\mu[\omega\wedge\omega']\geq \mu[\omega]\mu[\omega']\qquad\text{for all }\omega,\omega'\in\Omega_E,
\end{equation}
where $\omega\vee\omega'$ and $\omega\wedge\omega'$ refer to the coordinatewise maximum and minimum of $\omega$ and $\omega'$, respectively. 
If $\mu(\omega)>0$ for all $\omega\in\Omega_E$,~\eqref{eq:fkg_lattice} implies~\eqref{eq:fkg}~\cite{ForKasGin71,Hol74}.

In the setting of Section~\ref{sec:graph_repr_def}, we can identify $\gat_{\Lambda,K,K',K''}^{\#,\eta'}$ with $\gat_{\Lambda,K,K',K''}^{\#,\eta'}[\omega_E\in\cdot\,]$ and regard it as a measure on $\Omega_E$. With this identification, within a specific range of the coupling constants, the graphical representation satisfies the FKG lattice condition and, in particular, is positively associated.

\begin{lemma}\label{lem:gat_fkg}
Let $K,K',K''\in\R$ with $K\geq K''$ and $K>-K''$, and let $\Lambda\subset\Z^d$ finite. For $\#\in\{0,1\}$ and $\eta'\in\{+,\mrm{f}\}$, the measure $\gat_{\Lambda,K,K',K''}^{\#,\eta'}$ satisfies~\eqref{eq:fkg_lattice} if the weights~\eqref{eq:gat_weights_1} satisfy $\max\{\mrm{w}_2,\mrm{w}_3\}\leq 1$, which is equivalent to
\begin{equation}\label{eq:gat_fkg_restrictions}
K'\geq K''\quad\text{and}\quad \tanh K''\geq -\tanh K \tanh K'.
\end{equation}
\end{lemma}
These conditions imply $K'\geq 0$. It is useful to note that, once we additionally impose $K'\neq 0$, the assumptions of the above lemma are symmetric in $K$ and $K'$. Indeed, it is straightforward to check that they are then equivalent to 
\begin{equation}\label{eq:gat_fkg_regime}
\min\{K,K'\}>0,\quad \min\{K,K'\}\geq K''\quad\text{and}\quad \tanh K''\geq -\tanh K \tanh K'.
\end{equation}
As we will see below in Section~\ref{sec:atrc_coupling}, in the FKG regime $\mrm{w}_2,\mrm{w}_3\leq 1$, one may perform another FK-Ising-type expansion for the partition function in~\eqref{eq:gat_law_1}, which leads to a coupling of graphical representations for $\spin$ and $\spin'$. The resulting joint law is precisely that in~\cite{PfiVel97}. In this context, the statement of Lemma~\ref{lem:gat_fkg} is covered by~\cite[Proposition 4.1]{PfiVel97}, and we provide a proof in Appendix~\ref{sec:proof_fkg_lattice} for completeness.

\subsection{Coupling of pair of graphical representations in FKG regime}\label{sec:atrc_coupling}
In the setting of Section~\ref{sec:graph_repr_def}, fix $K,K',K''\in\R$ satisfying~\eqref{eq:gat_fkg_regime}, that is, the corresponding weights in~\eqref{eq:gat_weights_1} satisfy $\mrm{w}_1>0,\,\mrm{w}_2\in (0,1]$ and $\mrm{w}_3\in [0,1]$ with $(\mrm{w}_2,\mrm{w}_3)\neq (1,1)$. Let $\Lambda\subset\Z^d$ be finite, and set $E=\overline{\mbb{E}}_\Lambda$. Fix boundary conditions $\eta,\eta'\in\{+,\mrm{f}\}$, and set $\#=\ind{\eta=+}$ and $\#'=\ind{\eta'=+}$. Then there exist graphical representations for both configurations $s$ and $s'$. Performing an FK-Ising-type expansion for the partition function in~\eqref{eq:gat_law_1}, we construct a coupling of them. 

Define a measure on $\Omega_{\Z^d}\times\Omega_{\Z^d}$ by
\begin{equation}\label{eq:atrc_pair_law}
\atrc_{\Lambda,K,K',K''}^{\#,\#'}[\omega,\omega']\propto \mrm{u}_1^{\abs{\omega_E\setminus\omega'_E}}\,\mrm{u}_2^{\abs{\omega'_E\setminus\omega_E}}\,\mrm{u}_3^{\abs{\omega_E\cap\omega'_E}}\,2^{k_\Lambda(\omega)+k_\Lambda(\omega')}\,\ind{\omega\in\Omega_E^\#}\,\ind{\omega'\in\Omega_E^{\#'}},
\end{equation}
where the weights $\mrm{u}_i$ are given by
\begin{equation}\label{eq:atrc_weights_1}
\begin{aligned}
\mrm{u}_1:=\tfrac{\mrm{w}_1\mrm{w}_3}{\mrm{w}_2}=e^{2(K-K'')}-1,\qquad
\mrm{u}_2:=\tfrac{1}{\mrm{w}_2}-1=e^{2(K'-K'')}-1,\\
\mrm{u}_3:=\tfrac{\mrm{w}_1(1-\mrm{w}_3)}{\mrm{w}_2}=e^{2(K+K')}-e^{2(K-K'')}-e^{2(K'-K'')}+1.
\end{aligned}
\end{equation}

\begin{proposition}\label{prop:atrc_pair_coupling}
The measure $\atrc_{\Lambda,K,K',K''}^{\#,\#'}$ is a coupling of $\gat_{\Lambda,K,K',K''}^{\#,\eta'}$ and $\gat_{\Lambda,K',K,K''}^{\#',\eta}$. 
\end{proposition}
In Section~\ref{sec:sharp_af_iso}, we will see how this representation indicates certain monotonicity properties in the parameters of the isotropic model even when $K''=U$ is negative.
Observe that, in the case $K=K''$, one has $\mrm{u}_1=0$, whence $\omega_E\subseteq\omega'_E$ almost surely in the coupling~\eqref{eq:atrc_pair_law}, and analogously in the case $K'=K''$.
As announced before, the measure in~\eqref{eq:atrc_pair_law} is exactly that in~\cite{PfiVel97}, where it is also proved that the pair $(\omega,\omega')$ is jointly positively associated if $K''\geq 0$, while the same holds for $(\omega,-\omega')$ if $K''<0$.

\begin{proof}[Proof of Proposition~\ref{prop:atrc_pair_coupling}]
Consider the law in~\eqref{eq:gat_law_1} and the weights in~\eqref{eq:gat_weights_1}. Assume first that $\mrm{w}_3>0$, that is, $K>K''$. Writing $\abs{E_{\spin'}\setminus\omega}=\abs{E}-\abs{\omega_E}-\abs{E\setminus (\omega_E\cup E_{\spin'})}$ and $\abs{E_{\spin'}\cap\omega}=\abs{\omega_E}-\abs{\omega_E\setminus E_{\spin'}}$, we obtain
\begin{equation*}
\gat_{\Lambda,K,K',K''}^{\#,\eta'}[\omega]\propto\big(\tfrac{\mrm{w}_1\mrm{w}_3}{\mrm{w}_2}\big)^{\abs{\omega_E}}\,2^{k_\Lambda(\omega)}\,\ind{\omega\in\Omega_E^\#}\sum_{\spin'\in\Sigma_\Lambda^{\eta'}}\big(\tfrac{1}{\mrm{w}_2}\big)^{\abs{E\setminus (\omega_E\cup E_{\spin'})}}\big(\tfrac{1}{\mrm{w}_3}\big)^{\abs{\omega_E\setminus E_{\spin'}}}.
\end{equation*}
Now, writing $\mrm{w}_i^{-1}=1+(\mrm{w}_i^{-1}-1)$ for $i=2,3$, expanding and proceeding analogously to the proof of Proposition~\ref{prop:at_gat_laws}(ii), we deduce that
\begin{align*}
\sum_{\spin'\in\Sigma_\Lambda^{\eta'}}\big(\tfrac{1}{\mrm{w}_2}\big)^{\abs{E\setminus (\omega_E\cup E_{\spin'})}}\big(\tfrac{1}{\mrm{w}_3}\big)^{\abs{\omega_E\setminus E_{\spin'}}}&=\sum_{\substack{\omega'_1\subseteq E\setminus\omega_E\\ \omega'_2\subseteq\omega_E}}\big(\tfrac{1}{\mrm{w}_2}-1\big)^{\abs{\omega'_1}}\,\big(\tfrac{1}{\mrm{w}_3}-1\big)^{\abs{\omega'_2}}\,\sum_{\spin'\in\Sigma_\Lambda^{\eta'}}\ind{E_{\spin'}\cap(\omega'_1\cup\omega'_2)=\varnothing}\\
&\propto\sum_{\omega'\in\Omega_E^{\#'}}\big(\tfrac{1}{\mrm{w}_2}-1\big)^{\abs{\omega'_E\setminus\omega_E}}\,\big(\tfrac{1}{\mrm{w}_3}-1\big)^{\abs{\omega'_E\cap\omega_E}}\,2^{k_\Lambda(\omega')}.
\end{align*}
Therefore, $\gat_{\Lambda,K,K',K''}^{\#,\eta'}$ is the first marginal of $\atrc_{\Lambda,K,K',K''}^{\#,\#'}$. This remains valid when $\mrm{w}_3=0$ (that is, $K=K''$), which is shown by a similar calculation or by taking the limit $\mrm{w}_3\to 0$.
Then, by symmetry, the second marginal is given by $\gat_{\Lambda,K',K,K''}^{\#',\eta}$. 
\end{proof}

\section{Proof of Theorem~\ref{thm:antiferro}}\label{sec:af_phase}
The proof of Theorem~\ref{thm:antiferro} is a combination of classical energy-entropy arguments~\cite{Pei36} applied to the graphical representation (introduced in Section~\ref{sec:graph_repr_def}) of well chosen transformed variables.
We mention that, while the antiferromagnetic order in the product $\tau\tau'$ may be verified via a low temperature expansion and a Peierls argument without making use of the graphical representation, a high temperature expansion involves negative weights, see~\cite[Section 2.2]{PfiVel97}. Before proceeding to the proof, we introduce the \emph{spin-flip} maps and so-called \emph{Peierls contours}.
In regard to the latter, we follow~\cite[Chapter 5.7.4]{FriVel17} and make appropriate adjustments. 

\paragraph{Spin-flip map.}
For $A\subseteq\Z^d$, define $F_{A}:\Sigma_{\Z^d}\to\Sigma_{\Z^d}$ by setting
\begin{equation}\label{eq:sign_flip}
(F_{A}(\spin))_x=
\begin{cases}
		-\spin_x & \text{if }x\in A,\\
		\spin_x & \text{if }x\notin A.
\end{cases}
\end{equation}
For $A=\Z^d_{\mrm{odd}}$, we simply write $F_{\mrm{odd}}$.

\paragraph{Peierls contours.}
Associate to each $x\in \Z^d$ the closed unit hypercube $\mcal{S}_x:=x+[-\tfrac{1}{2},\frac{1}{2}]^d$.
Observe that the faces of $\mcal{S}_x,\,x\in\Z^d$, which are $(d-1)$-dimensional hypercubes, are in one--one correspondence with the edges of $\Z^d$.
We call these faces \emph{plaquettes}.
For $C\subset\Z^d$ finite and connected, the set $\mcal{M}(C):=\cup_{x\in C}\mcal{S}_x\subset\R^d$ is bounded and connected.
Therefore, there exists a unique connected component $\gamma(C)$ of the boundary $\partial\mcal{M}(C)$ (in the Euclidean sense) such that $C$ is contained in the bounded connected component of $\R^d\setminus\gamma(C)$.
Let $\mcal{E}(C)$ be the set of edges of $\Z^d$ that, when regarded as line-segments between their endvertices, intersect $\gamma(C)$.
Observe that, for a finite cluster $C$ of some $\omega\in\Omega_{\Z^d}$, the set $\mcal{E}(C)$ must be closed in $\omega$.
We say that a subset $F\subset\mbb{E}_{\Z^d}$ \emph{blocks} a vertex $x\in\Z^d$ if there exists a connected set $C\subset\Z^d$ with $x\in C$ and $\mcal{E}(C)=F$. We call $F\subset\mbb{E}_{\Z^d}$ \emph{blocking} if it blocks some vertex $x\in\Z^d$.

In the proof of Theorem~\ref{thm:antiferro}, we will need an upper bound on the number of blocking $F\subset E$ with $\abs{F}=k$ for some fixed $E\subset\mbb{E}_{\Z^d}$.
We endow $E$ with a graph structure by declaring two edges in $E$ to be connected if the corresponding plaquettes share a $(d-2)$-dimensional hypercube.
It is easy to see that the resulting graph has maximum degree $6(d-1)$.
By construction, the set of blocking $F\subset E$ of size $k$ is embedded in the set of connected subgraphs with $k$ vertices of this graph. Using for example~\cite[Lemma 3.38]{FriVel17}, it is easy to derive the following crude bound.
\begin{lemma}\label{lem:bound_block}
Let $E\subseteq\mbb{E}_{\Z^d}$ finite and $k\geq 1$. The number of blocking $F\subseteq E$ with $\abs{F}=k$ is bounded by
\[
\abs{E}(6(d-1))^{6(d-1)k}.
\] 
\end{lemma}
We are now ready to prove the theorem.

\begin{proof}[Proof of Theorem~\ref{thm:antiferro}]
Fix $d\geq 2$ and $n\geq 1$, and set $E=\overline{\mbb{E}}_{\msf{B}_n}$.
Assume that $-U>J>0$ and let $(\tau,\tau')$ be distributed according to $\at_{\msf{B}_n,J,U}^{+,\pm}$. Observe that the statement of the theorem is trivial for $x\in\Z^d\setminus\msf{B}_n$, whence we fix $x\in\msf{B}_n$.
The proof is divided into two steps.

\noindent\textbf{Step 1.} For $J\geq 0$ small enough, there exists $c=c(d,J)>0$ such that $\langle\tau_x\rangle_{\msf{B}_n,J,U}^{+,\pm}\leq e^{-c(n-\norm{x})}$.

First, we apply the change of variables $\spin=\tau,\,\spin'=\tau\tau'$, and observe that $(\spin,\spin')$ is distributed according to $\at_{\msf{B}_n,J,U,J}^{+,\pm}$. The consequence, Corollary~\ref{cor:gat_corr_conn}, of Proposition~\ref{prop:at_gat_laws}(i), applied to $(\spin,\spin')$ and with $K=K''=J,\,K'=U$, implies
\[
\langle \tau_x\rangle_{\msf{B}_n,J,U}^{+,\pm}=\gat_{\msf{B}_n,J,U,J}^{1,\pm}[x\leftrightarrow \Z^d\setminus\msf{B}_n]\leq \sum_{\Gamma:\,x\to\Z^d\setminus\msf{B}_n}\gat_{\msf{B}_n,J,U,J}^{1,\pm}[\Gamma\text{ open}],
\]
where we used a union bound, and where the summation index $\Gamma:x\to\Z^d\setminus\msf{B}_n$ refers to all paths connecting $x$ to $\Z^d\setminus\msf{B}_n$ in $\overline{\msf{B}_n}$. Recall Lemma~\ref{lem:finite_energy}, and observe that, by our assumption on $J,U$, the weights in~\eqref{eq:gat_weights_1} satisfy $\mrm{w}_2>1$ and $\mrm{w}_3=0$. Hence, for a path $\Gamma=(x_0,\dots,x_\ell)$ with $e_i=x_ix_{i+1}$, Lemma~\ref{lem:finite_energy} gives
\begin{equation}\label{eq:atrc_path_bound_1}
\gat_{\msf{B}_n,J,U,J}^{1,\pm}[\Gamma\text{ open}]=\prod_{i=0}^{\ell-1}\gat_{\msf{B}_n,J,U,J}^{1,\pm}[\omega_{e_i}=1\,|\,\omega_{e_j}=1\,\forall\,j<i]\leq \big(1+\mrm{w}_1^{-1}\big)^{-\ell}.
\end{equation}
Bounding the number of $\Gamma:x\to\Z^d\setminus\msf{B}_n$ of length $\ell$ crudely by $(2d)^\ell$, we obtain
\[
\langle \tau_0\rangle_{\msf{B}_n,J,U}^{+,\pm}\leq\sum_{\ell>n-\norm{x}}\,\sum_{\substack{\Gamma:\,x\to\Z^d\setminus\msf{B}_n,\\ \abs{\Gamma}=\ell}}\gat_{\msf{B}_n,J,U,J}^{1,\pm}[\Gamma\text{ open}]\leq \sum_{\ell>n-\norm{x}} \Big(\tfrac{2d}{1+\mrm{w}_1^{-1}}\Big)^\ell,
\]
and since $\mrm{w}_1=e^{4J}-1$ tends to zero as $J\to 0$, the claim follows.

\noindent\textbf{Step 2.} For $-J-U$ large enough, there exists $c=c(d,J+U)>0$ such that
\[
\langle\tau_x\tau'_x\rangle_{\msf{B}_n,J,U}^{+,\pm}
\begin{cases}
\geq c & \text{if }x\in\Z^d_{\mrm{even}},\\
\leq -c & \text{if }x\in\Z^d_{\mrm{odd}},
\end{cases}
\]
and $\mcal{C}$ with the desired properties exists.

Recall the definition of $F_{\mrm{odd}}$ in~\eqref{eq:sign_flip}, and apply the change of variables $\spin=F_{\mrm{odd}}(\tau\tau')$, $\spin'=\tau$.
Observe that the law of $(\spin,\spin')$ is given by $\at_{\msf{B}_n,-U,J,-J}^{+,+}$ and that
\[
\langle\tau_x\tau'_x\rangle_{\msf{B}_n,J,U}^{+,\pm}=\big(\ind{x\in\Z^d_\mrm{even}}-\ind{x\in\Z^d_\mrm{odd}}\big)\langle\spin_x\rangle_{\msf{B}_n,-U,J,-J}^{+,+}.
\]
Let $\widetilde{\omega}$ be distributed according to $\gat_{\msf{B}_n,-U,J,-J}^{1,+}$, and denote by $\widetilde{\mrm{w}}_i$ the corresponding weights in~\eqref{eq:gat_weights_1} with $K=-U,\,K'=-K''=J$. Define $\mcal{C}_x=\mcal{C}_x(\widetilde{\omega})\subseteq\Z^d$ to be the cluster of $x$ in $\widetilde{\omega}$. Then, by Corollary~\ref{cor:gat_corr_conn},
\[
\langle\spin_x\rangle_{\msf{B}_n,-U,J,-J}^{+,+}=1-\gat_{\msf{B}_n,-U,J,-J}^{1,+}[\mcal{E}(\mcal{C}_x)\subseteq E]\geq 1-\sum_{F\subseteq E}\gat_{\msf{B}_n,-U,J,-J}^{1,+}[F\text{ closed}],
\]
where we used that $\mcal{C}_x\subseteq\msf{B}_n$ if and only if $\mcal{E}(\mcal{C}_x)\subseteq E$, and where the summation is over all non-empty $F\subseteq E$ that block~$x$. The assumption $-U>J>0$ implies that $\widetilde{\mrm{w}}_2<1<\widetilde{\mrm{w}}_3$.
Analogously to Step 1, for any $F\subseteq E$, Lemma~\ref{lem:finite_energy} gives
\begin{equation}\label{eq:bound_gat_edge_bdry_closed}
\gat_{\msf{B}_n,-U,J,-J}^{1,+}[F\text{ closed}]\leq \big(1+\tfrac{\widetilde{\mrm{w}}_1}{2}\big)^{-\abs{F}}.
\end{equation}
Furthermore, any $F\subseteq E$ with $\abs{F}=k$ that blocks $x$ necessarily satisfies $F\subseteq\mbb{E}_{\msf{B}_k(x)}=:E^x_k$.
Altogether, we obtain
\[
\langle\spin_x\rangle_{\msf{B}_n,-U,J,-J}^{+,+}\geq 1-\sum_{k\geq 2d}\sum_{\substack{F\subseteq E:\\ \abs{F}=k}}\big(1+\tfrac{\widetilde{\mrm{w}}_1}{2}\big)^{-k}\geq 1-\sum_{k\geq 2d}\abs{E^x_k}(6(d-1))^{6(d-1)k}\big(1+\tfrac{\widetilde{\mrm{w}}_1}{2}\big)^{-k},
\]
where we applied Lemma~\ref{lem:bound_block} to bound the number of blocking $F\subseteq E^x_k$ with $\abs{F}=k$. Since $\abs{E^x_k}=\mcal{O}(k^d)$ and $\widetilde{\mrm{w}}_1=e^{-2(J+U)}-1\to\infty$ as $J+U\to-\infty$, the first claim of Step 2 is proved.

Finally, define $\mcal{C}\subset\Z^d$ as the connected component of $\Z^d\setminus\msf{B}_n$ in $\widetilde{\omega}\sim\gat_{\msf{B}_n,-U,J,-J}^{1,+}$, that is, 
\[
\mcal{C}:=\{y\in\Z^d: y\xleftrightarrow{\widetilde{\omega}}\Z^d\setminus\msf{B}_n\}.
\]
Then, by Proposition~\ref{prop:at_gat_laws}(i), $\spin=F_\mrm{odd}(\tau\tau')$ is constant $+1$ on $\mcal{C}$, and hence $\tau\tau'$ is constant $+1$ on $\mcal{C}_\mrm{even}$ and $-1$ on $\mcal{C}_\mrm{odd}$. 
It remains to show~\eqref{eq:at_af_log_sized_holes} for $c=c(J+U)$ when $-J-U$ is large enough. If $\Lambda\subseteq\msf{B}_n\setminus\mcal{C}$ is connected, we can assume without loss of generality that $\mcal{E}(\Lambda)\cap\widetilde{\omega}=\varnothing$ (otherwise, $\Lambda$ may be enlarged).
Noticing that $\abs{\mcal{E}(\Lambda)}$ is greater than $\mrm{diam}(\Lambda)$, bounding the number of blocking $F\subseteq E$ of a given size by applying Lemma~\ref{lem:bound_block}, and using the bound~\eqref{eq:bound_gat_edge_bdry_closed}, we deduce
\[
\at_{\msf{B}_n,J,U}^{+,\pm}[\exists\,\Lambda\subseteq\msf{B}_n\setminus\mcal{C}\text{ connected with }\mrm{diam}(\Lambda)> \log n]\leq\sum_{k>\log n}\abs{E}(6(d-1))^{6(d-1)k}\big(1+\tfrac{\widetilde{\mrm{w}}_1}{2}\big)^{-k}. 
\]
Since $\abs{E}=\mcal{O}(n^d)$ and $\widetilde{\mrm{w}}_1\to\infty$ as $J+U\to-\infty$, the proof is complete.
\end{proof}

\section{Coupling with the eight-vertex model}\label{sec:mappings_models}
Duality transformations applied to pairs of AT spin configurations relate the AT model on the square lattice $\Z^2$ and on its dual $(\Z^2)^*$~\cite{MitSte71,Fan72c,PfiVel97}.
On the other hand, keeping one of the AT configurations and applying such a transformation to one of the other, results in a pair of spin configurations, one of which is defined on $\Z^2$ and the other on $(\Z^2)^*$.
These configurations then follow the law of the spin representation~\cite{Wu71,KadWeg71,GlaMan21,Lis22,GlaPel23} of a staggered (i.e. with non-translation-invariant weights~\cite{WuLin75,HsuLinWu75}) eight-vertex model~\cite{Sut70,FanWu70} on the medial graph.
In special cases, the latter may be non-staggered or reduces to the six-vertex model~\cite{Pau35,Rys63}.

The relation between the AT and the eight-vertex models on the square lattice has first been noticed in~\cite{Fan72} comparing their critical properties, and it was made explicit in~\cite{Fan72b,Weg72} on a level of partition functions (see also~\cite{Wu77,HuaDenJacSal13} and~\cite[Chapter 12.9]{Bax89}).
In~\cite{GlaPel23}, a coupling of the self-dual isotropic AT model and the six-vertex model (more precisely, the F-model~\cite{Rys63}) was constructed.

In Section~\ref{sec:at_8v_coupling}, we use the graphical representation introduced in Section~\ref{sec:graph_repr_def} to unify the above approaches and construct couplings of the AT and eight-vertex measures.
In Section~\ref{sec:8v_ferroel_phase}, we explain the relation between Theorems~\ref{thm:antiferro}~and~\ref{thm:hf_loc} and we prove the latter. 

\subsection{Duality coupling}\label{sec:at_8v_coupling}
\begin{figure}
\begin{center}
\includegraphics[page=1,height=7cm]{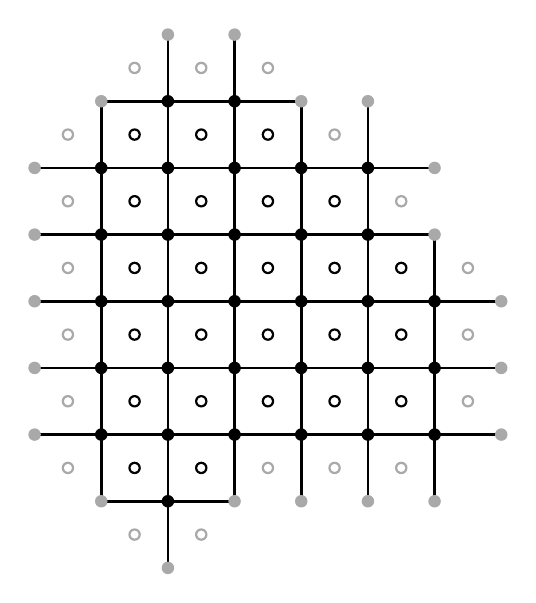}
\hspace{25pt}
\includegraphics[page=2,height=7cm]{square_lattice_domain.pdf}
\end{center}
\caption{\textit{Left:} A domain $\Lambda\subset\Z^2$ (black solid vertices) with external boundary $\overline{\Lambda}\setminus\Lambda$ (grey solid vertices), the edge set $E:=\overline{\mbb{E}}_\Lambda$ (black solid edges), and the dual $\Lambda^*\subset(\Z^2)^*$ (black hollow vertices) with external boundary $\overline{\Lambda^*}\setminus\Lambda^*$ (grey hollow vertices).
\textit{Right:} An edge configuration $\omega\in\Omega_E^0$ (black solid edges) and the relevant part of its dual $\omega^*$ (black and grey dashed edges, the grey ones are dual to edges outside $E$).}
\label{fig:domain_dual}
\end{figure}

Before constructing the coupling, we introduce some notation concerning the dual square lattice $(\Z^2)^*$ and discuss the duality map for percolation configurations.
\vspace{-5pt}
\paragraph{Dual square lattice.}
The dual graph of the square lattice has vertex-set $(\Z^2)^*=(\tfrac{1}{2},\tfrac{1}{2})+\Z^2$ and edges between nearest neighbours, and we again write $(\Z^2)^*$ for the graph itself. 
Regarded as line-segments between their endvertices, each edge $e\in\mbb{E}_{\Z^2}$ intersects a unique edge $e^*\in\mbb{E}_{(\Z^2)^*}$, and we call $e^*$ the edge dual to $e$ and vice versa. For $E\subseteq\mbb{E}_{\Z^2}$ and $\spin\in\Sigma_{V'}$ with $V'\subseteq(\Z^2)^*$, write $E_{\spin}$ for the set of edges in $E$ dual to disagreement edges of $\spin$, that is, 
\begin{equation*}
E_{\spin}:=\{e\in E: e^*=x'y'\text{ and }\spin_{x'}\neq \spin_{y'}\}.
\end{equation*}
With each percolation configuration $\omega\in\Omega_{\Z^2}$ is associated a dual percolation configuration $\omega^*\in\Omega_{(\Z^2)^*}$ given by $\omega^*_{e^*}=1-\omega_e,\,e\in\mbb{E}_{\Z^2}$.
We call a finite subset $\Lambda\subset\Z^2$ a \emph{domain} if all bounded faces of the graph $(\overline{\Lambda},\overline{\mbb{E}}_\Lambda)$ are unit squares, as illustrated in Figure~\ref{fig:domain_dual}.
In this case, we define its dual $\Lambda^*\subset(\Z^2)^*$ as the set of vertices in $(\Z^2)^*$ at which the bounded faces of $(\overline{\Lambda},\overline{\mbb{E}}_\Lambda)$ are centred, as depicted in the figure. 
Then, for $E=\overline{\mbb{E}}_\Lambda$ and $\omega\in\Omega_E^0$, the set of faces of $(\overline{\Lambda},\omega)$ is in one--one correspondence with the set of clusters of $\omega^*$ that intersect $\overline{\Lambda^*}$ (see Figure~\ref{fig:domain_dual}), and Euler's formula for planar graphs implies
\begin{equation}\label{eq:eulers_formula_consequence}
k_\Lambda(\omega)=\abs{\overline{\Lambda}}-\abs{\omega_E}+k_{\Lambda^*}(\omega^*)-1,
\end{equation}
see, e.g.,~\cite[Chapter 6.1]{Gri06}.

\begin{figure}
    \centering
    \includegraphics[scale=0.6]{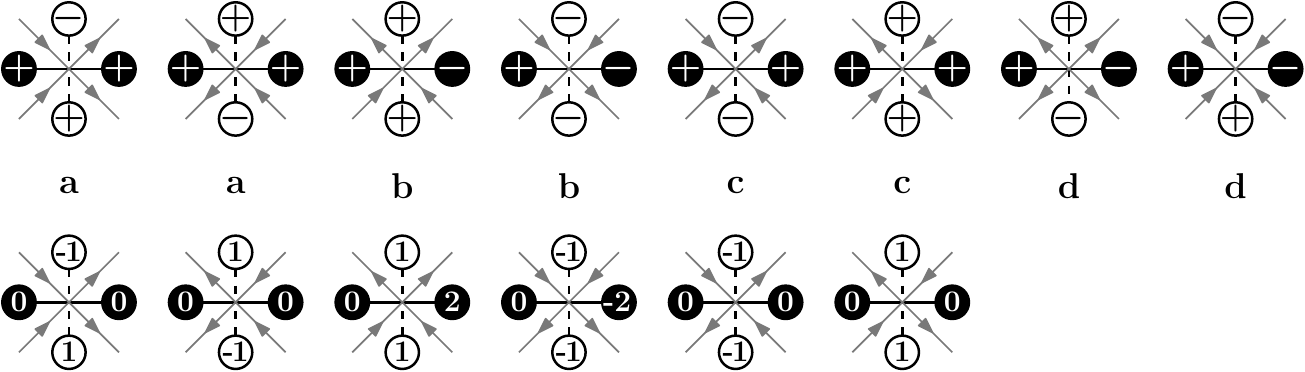}
    \caption{\textit{Top:} the eight possible local eight-vertex spin configurations at an horizontal edge $e$ of $\Z^2$ (solid black) and its dual $e^*$ (dashed black) when the spin $\sigma^\bullet$ is fixed to be $+1$ at the left endpoint of $e$. The corresponding orientations of the edges of the medial graph (solid gray) are obtained by orienting an edge so that the vertex of $\Z^2$ (black) is on its left precisely if the spins separated by the edge coincide.
    \textit{Centre:} the corresponding local weights. For the same sequence of edge orientations at a vertical edge of $\Z^2$, the weights $\mbf{a}$ and $\mbf{b}$ are interchanged. 
    \textit{Bottom:} the height function values corresponding to types 1-6 obtained by~\eqref{eq:hf_gradient} when the height is fixed to be $0$ at the left endpoint of $e$. In these cases, the edge orientations are visualisations of the discrete gradient of the height function.}
    \label{fig:eight_vertex_types}
\end{figure}

\paragraph{Duality coupling.}
We formulate the coupling in the setting of Section~\ref{sec:graph_repr_def}. Let $K,K',K''\in\R$ with $K\geq\abs{K''}$, and define the corresponding eight-vertex weights
\begin{equation}\label{eq:8v_weights}
\begin{aligned}
\mbf{a}&=e^{2(K+K'')}-1,\qquad		&		\mbf{b}&=e^{-2K'}\big(e^{2K}+e^{2K''}\big),\\
\quad \mbf{c}&=e^{2(K+K'')}+1,\qquad			&		\mbf{d}&=e^{-2K'}\big(e^{2K}-e^{2K''}\big).
\end{aligned}
\end{equation}
Let $\Lambda\subset\Z^2$ be a domain, and set $E=\overline{\mbb{E}}_{\Lambda}$.
For convenience regarding Theorem~\ref{thm:hf_loc}, we take free boundary conditions $\eta=\mrm{f}$ and arbitrary $\eta'\in\Sigma_{\Z^2}$. 
Take $(\spin,\spin',\omega)$ distributed as in Proposition~\ref{prop:at_gat_laws}(i). Sample a pair of spin configurations $\sigma=(\sigma^\bullet,\sigma^\circ)\in\Sigma_{\Lambda}^{\eta'}\times\Sigma_{\Lambda^*}^+$ as follows:
\begin{itemize}
\setlength\itemsep{0.1em}
\item set $\sigma^\bullet=\spin'$, 
\item sample $\sigma^\circ$ by assigning $+1$ to the cluster of $\omega^*$ that intersects $(\Z^2)^*\setminus\Lambda^*$ and $\pm1$ uniformly and independently to the clusters of $\omega^*$ that are contained in $\Lambda^*$.
\end{itemize}
We denote the law of $\sigma$ on $\Sigma_{\Z^2}\times\Sigma_{(\Z^2)^*}$ by $\eightv_{\Lambda,\mbf{a},\mbf{b},\mbf{c},\mbf{d}}^{\eta',+}$, and call it the staggered eight-vertex (spin) measure on $\Lambda\cup\Lambda^*$ with parameters $\mbf{a},\mbf{b},\mbf{c},\mbf{d}$ and boundary condition $(\eta',+)$.
Recall that we conventionally use the same symbols for random variables and their deterministic realisations.
\begin{proposition}\label{prop:at_atrc_8v_coupling}
\noindent(i) The joint law of $(\spin,\spin',\omega,\sigma)$ on $(\Sigma_{\Z^2})^2\times\Omega_{\Z^2}\times(\Sigma_{\Z^2}\times\Sigma_{(\Z^2)^*})$ is given by
\begin{equation}\label{eq:at_atrc_8v_joint_law}
\begin{aligned}
\Prob{\spin,\spin',\omega,\sigma}\propto e^{2(K''-K')\abs{E_{\spin'}}}\,\big(e^{2(K-K'')}-1\big)^{\abs{\omega_E\cap E_{\spin'}}}\,\big(e^{2(K+K'')}-1\big)^{\abs{\omega_E\setminus E_{\spin'}}}\,2^{-k_{\Lambda^*}(\omega^*)}\\
\cdot\,\ind{(\spin,\spin')\in\Sigma_\Lambda^{\mrm{f}}\times\Sigma_\Lambda^{\eta'}}\,\ind{\omega\in\Omega_E^0}\,\ind{\omega\cap E_\spin=\varnothing}\,\ind{\sigma^\bullet=\spin'}\,\ind{\sigma^\circ\in\Sigma_{\Lambda^*}^+}\,\ind{E_{\sigma^\circ}\subseteq\omega}.
\end{aligned}
\end{equation}
It is a coupling of
\[(\spin,\spin')\sim\at_{\Lambda,K,K',K''}^{\mrm{f},\eta'},\quad\omega\sim\gat_{\Lambda,K,K',K''}^{0,\eta'}\quad\text{and}\quad\sigma\sim\eightv_{\Lambda,\mbf{a},\mbf{b},\mbf{c},\mbf{d}}^{\eta',+},
\]
where the relation between $K,K',K''$ and $\mbf{a},\mbf{b},\mbf{c},\mbf{d}$ is given by~\eqref{eq:8v_weights}.

\noindent(ii) The law of $\sigma$ on $\Sigma_{\Z^2}\times\Sigma_{(\Z^2)^*}$ is given by
\begin{equation}\label{eq:8v_law}
\eightv_{\Lambda,\mbf{a},\mbf{b},\mbf{c},\mbf{d}}^{\eta',+}[\sigma]\propto
\mbf{a}^{\abs{E_{\sigma^\circ}\setminus E_{\sigma^\bullet}}}\,
\mbf{b}^{\abs{E_{\sigma^\bullet}\setminus E_{\sigma^\circ}}}\,
\mbf{c}^{\abs{E\setminus (E_{\sigma^\bullet}\cup E_{\sigma^\circ})}}\,
\mbf{d}^{\abs{E_{\sigma^\bullet}\cap E_{\sigma^\circ}}}\,\ind{\sigma^\bullet\in\Sigma_\Lambda^{\eta'}}\,\ind{\sigma^\circ\in\Sigma_{\Lambda^*}^+}.
\end{equation}
\end{proposition}
To obtain a relation with the classical representation of the eight-vertex model in terms of edge orientations, one may apply the maps in~\cite{Wu71,KadWeg71} or~\cite{Lis22,GlaPel23}. In the above setting, these are orientations of the edges of the medial graph of $\Z^2$, which is a rotated square lattice. We will only comment on this representation in Figure~\ref{fig:eight_vertex_types}, and we focus on those in terms of spin configurations and height functions.
We remark that the model is non-staggered, meaning in our context that the same weights are assigned to disagreements in $\sigma^\bullet$ and $\sigma^\circ$, precisely if $\mbf{a}=\mbf{b}$, which corresponds to the self-dual manifold of the AT model~\cite{MitSte71,Fan72c}.

We emphasise that our choice of weights differs from the classical one in~\cite{Bax89}. In fact, the weights $(a,b,c,d)$ in~\cite[Equation (10.2.1)]{Bax89} coincide with our weights $(\mbf{a},\mbf{b},\mbf{c},\mbf{d})$ on one sublattice of the medial graph and with $(\mbf{b},\mbf{a},\mbf{c},\mbf{d})$ on the other, see Figure~\ref{fig:eight_vertex_types}.

\begin{proof}[Proof of Proposition~\ref{prop:at_atrc_8v_coupling}]
\noindent(i) While $\sigma^\bullet$ is uniquely determined by $\spin'$, there are $2^{k_{\Lambda^*}(\omega^*)-1}$ possibilities for $\sigma^\circ\in\Sigma_{\Lambda^*}^+$ with $E_{\sigma^\circ}\subseteq\omega$ each of which is equally likely.
Hence, the joint law of $(\spin,\spin',\omega,\sigma)$ is given by~\eqref{eq:at_atrc_8v_joint_law}.

\noindent(ii) We first sum~\eqref{eq:at_atrc_8v_joint_law} over $\spin\in\Sigma_\Lambda^{\mrm{f}}$ with $\omega\cap E_\spin=\varnothing$, of which there are $2^{k_\Lambda(\omega)}$, and then over $\spin'\in\Sigma_\Lambda^{\eta'}$, which is uniquely determined by $\spin'=\sigma^\bullet$, to obtain
\begin{equation*}
\begin{aligned}
\Prob{\omega,\sigma}&\propto e^{2(K''-K')\abs{E_{\sigma^\bullet}}}\,\big(e^{2(K-K'')}-1\big)^{\abs{\omega_E\cap E_{\sigma^\bullet}}}\,\big(e^{2(K+K'')}-1\big)^{\abs{\omega_E\setminus E_{\sigma^\bullet}}}\,2^{k_{\Lambda}(\omega)-k_{\Lambda^*}(\omega^*)}\\
&\qquad\cdot\ind{\omega\in\Omega_E^0}\,\ind{\sigma^\bullet\in\Sigma_\Lambda^{\eta'}}\,\ind{\sigma^\circ\in\Sigma_{\Lambda^*}^+}\,\ind{E_{\sigma^\circ}\subseteq\omega}\\
&\propto e^{2(K''-K')\abs{E_{\sigma^\bullet}}}\,\Big(\tfrac{e^{2(K-K'')}-1}{2}\Big)^{\abs{\omega_E\cap E_{\sigma^\bullet}}}\,\Big(\tfrac{e^{2(K+K'')}-1}{2}\Big)^{\abs{\omega_E\setminus E_{\sigma^\bullet}}}\\
&\qquad\cdot\ind{\omega\in\Omega_E^0}\,\ind{\sigma^\bullet\in\Sigma_\Lambda^{\eta'}}\,\ind{\sigma^\circ\in\Sigma_{\Lambda^*}^+}\,\ind{E_{\sigma^\circ}\subseteq\omega},
\end{aligned}
\end{equation*}
where we also used that $k_\Lambda(\omega)-k_{\Lambda^*}(\omega^*)=-\abs{\omega_E}+\mrm{const}(\Lambda)$, see~\eqref{eq:eulers_formula_consequence} and above.
Summing over $\omega=E_{\sigma^\circ}\cup\omega_1\cup\omega_2$ with $\omega_1\subseteq E_{\sigma^\bullet}\setminus E_{\sigma^\circ}$ and $\omega_2\cap (E_{\sigma^\bullet}\cup E_{\sigma^\circ})=\varnothing$,
\begin{align*}
\Prob{\sigma}&\propto e^{2(K''-K')\abs{E_{\sigma^\bullet}}}\,\Big(\tfrac{e^{2(K-K'')}-1}{2}\Big)^{\abs{E_{\sigma^\circ}\cap E_{\sigma^\bullet}}}\,\Big(\tfrac{e^{2(K+K'')}-1}{2}\Big)^{\abs{E_{\sigma^\circ}\setminus E_{\sigma^\bullet}}}\,\ind{\sigma^\bullet\in\Sigma_\Lambda^{\eta'}}\,\ind{\sigma^\circ\in\Sigma_{\Lambda^*}^+}\\
&\qquad\cdot\sum_{\omega_1\subseteq E_{\sigma^\bullet}\setminus E_{\sigma^\circ}}\Big(\tfrac{e^{2(K-K'')}-1}{2}\Big)^{\abs{\omega_1}}\,\sum_{\omega_2\subseteq E\setminus (E_{\sigma^\bullet}\cup E_{\sigma^\circ})}\Big(\tfrac{e^{2(K+K'')}-1}{2}\Big)^{\abs{\omega_2}}\\
&= e^{2(K''-K')\abs{E_{\sigma^\bullet}}}\,\Big(\tfrac{e^{2(K-K'')}-1}{2}\Big)^{\abs{E_{\sigma^\circ}\cap E_{\sigma^\bullet}}}\,\Big(\tfrac{e^{2(K+K'')}-1}{2}\Big)^{\abs{E_{\sigma^\circ}\setminus E_{\sigma^\bullet}}}\,\ind{\sigma^\bullet\in\Sigma_\Lambda^{\eta'}}\,\ind{\sigma^\circ\in\Sigma_{\Lambda^*}^+}\\
&\qquad\cdot \Big(\tfrac{e^{2(K-K'')}+1}{2}\Big)^{\abs{E_{\sigma^\bullet}\setminus E_{\sigma^\circ}}}\,\Big(\tfrac{e^{2(K+K'')}+1}{2}\Big)^{\abs{E\setminus (E_{\sigma^\bullet}\cup E_{\sigma^\circ})}}.
\end{align*}
Collecting like terms and multiplying by $2^{\abs{E}}$, we arrive at~\eqref{eq:8v_law}, and the proof is complete.
\end{proof}
The six-vertex model can be considered a special case of the eight-vertex model. Its characterisation is based on the \emph{ice-rule}, which guarantees that the \emph{height function} observable is well-defined.
Via the map~\cite{Lis22,GlaPel23}, the ice-rule states that $E_{\sigma^\bullet}$ and $E_{\sigma^\circ}$ must be disjoint, which is to say that
\begin{equation}\label{eq:ice_rule}
(\sigma^\bullet_x-\sigma^\bullet_y)(\sigma^\circ_{x'}-\sigma^\circ_{y'})=0\qquad\text{for all }e=xy\in\mbb{E}_{\Z^2}\text{ with }e^*=x'y'.
\end{equation}
In the above notation, this is equivalent to $\mbf{d}=0$, which in turn is equivalent to $K=K''$. We then write $\sixv_{\Lambda,\mbf{a},\mbf{b},\mbf{c}}^{\eta',+}$ for the measure $\eightv_{\Lambda,\mbf{a},\mbf{b},\mbf{c},0}^{\eta',+}$, and we call it the staggered six-vertex (spin) measure on $\Lambda\cup\Lambda^*$ with parameters $\mbf{a},\mbf{b},\mbf{c}$ and boundary condition $(\eta',+)$.

\paragraph{Height function representation of the six-vertex model.}
Recall the definitions of a quad and a (six-vertex) height function given in Section~\ref{sec:def_models}.
A configuration $\sigma=(\sigma^\bullet,\sigma^\circ)\in\Sigma_{\Z^2}\times\Sigma_{(\Z^2)^*}$ that satisfies the ice-rule~\eqref{eq:ice_rule} defines a height function $h$ up to an additive constant: for $x\in\Z^2$ and $x'\in(\Z^2)^*$ belonging to the same quad, define
\begin{equation}\label{eq:hf_gradient}
h_{x'}-h_x=\sigma^\bullet_x\,\sigma^\circ_{x'}.
\end{equation}
The ice-rule~\eqref{eq:ice_rule} guarantees that the sum of the four gradients~\eqref{eq:hf_gradient} around each quad is zero.
In particular, there is a one--one correspondence between such configurations $\sigma$ and height functions that take a fixed value at some fixed vertex of $\Z^2$ or $(\Z^2)^*$, see Figure~\ref{fig:eight_vertex_types}.
Given $h$ and $\sigma$ related by~\eqref{eq:hf_gradient}, observe that $h$ is constant on a fixed edge of $\Z^2$ or $(\Z^2)^*$ if and only if $\sigma$ is:
\begin{equation}\label{eq:hf_spin_agree}
\begin{aligned}
\text{for any }xy\in\mbb{E}_{\Z^2},\quad h_x=h_y\quad&\text{if and only if}\quad \sigma^\bullet_x=\sigma^\bullet_y,\\
\text{for any }x'y'\in\mbb{E}_{(\Z^2)^*},\quad h_{x'}=h_{y'}\quad&\text{if and only if}\quad \sigma^\circ_{x'}=\sigma^\circ_{y'}.
\end{aligned}
\end{equation}

\subsection{Relation of Theorems~\ref{thm:antiferro} and~\ref{thm:hf_loc} and proof of Theorem~\ref{thm:hf_loc}}
\label{sec:8v_ferroel_phase}

Recall the definition of the $\pm$-boundary condition given in Section~\ref{sec:results}.
In order to understand the relationship between Theorems~\ref{thm:antiferro}~and~\ref{thm:hf_loc}, consider the coupling~\eqref{eq:at_atrc_8v_joint_law} with $\spin=\tau,\,\spin'=\tau\tau'$, where $\tau$ and $\tau\tau'$ are as in Theorem~\ref{thm:antiferro} but (for convenience) with boundary conditions taken free for $\tau$ and $\pm$ for $\tau\tau'$. The law of the corresponding $\sigma$ is then given by $\sixv_{\msf{B}_n,\mbf{a},\mbf{b},\mbf{c}}^{\pm,+}$. Since $\sigma^\bullet=\tau\tau'$ is antiferromagnetically ordered, the corresponding even heights on $\Z^2$ vary along edges of $\Z^2$. On the other hand, $\tau$ is disordered, admitting exponential decay of correlations, whence its dual $\sigma^\circ$ and therefore the corresponding odd heights on $(\Z^2)^*$ are ordered. The latter is sufficient to keep the height function flat, meaning its variance is uniformly bounded.
The first step is to identify the law of the corresponding height function pinned to $1$ on $(\Z^2)^*\setminus\msf{B}_n^*$ as the measure $\hf_{\msf{B}_n,\mbf{a},\mbf{b},\mbf{c}}^{0/2,1}$ from Theorem~\ref{thm:hf_loc}.

\begin{lemma}\label{lem:hf_law}
Let $\mbf{a},\mbf{b},\mbf{c}>0$ and $n\geq 1$. Let $\sigma=(\sigma^\bullet,\sigma^\circ)$ be distributed according to $\sixv_{\msf{B}_n,\mbf{a},\mbf{b},\mbf{c}}^{\pm,+}$. There exists a unique random height function~$h$ related to $\sigma$ by~\eqref{eq:hf_gradient} and satisfying~$h=1$ on $(\Z^2)^*\setminus\msf{B}_n^*$. Its law is given by $\hf_{\msf{B}_n,\mbf{a},\mbf{b},\mbf{c}}^{0/2,1}$.
\end{lemma}

\begin{proof}
Since $\sigma^\circ\in\Sigma_{\msf{B}_n^*}^+$, there is a unique height function $h$ satisfying~\eqref{eq:hf_gradient} and $h=1$ on $(\Z^2)^*\setminus\msf{B}_n^*$. Clearly, equation~\eqref{eq:hf_gradient} and $\sigma^\bullet\in\Sigma_{\msf{B}_n}^\pm$ then imply that $h=0$ on $\Z^2_\mrm{even}\setminus\msf{B}_n$ and $h=2$ on $\Z^2_\mrm{odd}\setminus\msf{B}_n$.

Let $E=\overline{\mbb{E}}_{\msf{B}_n}$. By~\eqref{eq:hf_spin_agree}, we have $E_{\sigma^\bullet}=E_{h^\bullet}$ and $E_{\sigma^\circ}=E_{h^\circ}$. Altogether, we conclude that the law of $h$ is given by~\eqref{eq:hf_meas_def} with $t=2\mathbbm{1}_{\Z^2_{\mrm{odd}}}+\mathbbm{1}_{(\Z^2)^*}$, and the proof is complete. 
\end{proof}

Although some parts of the proof of Theorem~\ref{thm:hf_loc} are now analogous to that of Theorem~\ref{thm:antiferro}, it is not a direct consequence, whence we provide it for completeness.

\begin{proof}[Proof of Theorem~\ref{thm:hf_loc}]
Fix $\mbf{a},\mbf{b},\mbf{c}>0$ and $n\geq 1$, and set $E=\overline{\mbb{E}}_{\msf{B}_n}$. Consider the coupling of
\[
(\tau,\tau\tau')\sim\at_{\msf{B}_n,J,U,J}^{\mrm{f},\pm},\quad\omega\sim\gat_{\msf{B}_n,J,U,J}^{0,\pm}\quad\text{and}\quad\sigma=(\sigma^\bullet,\sigma^\circ)\sim\sixv_{\msf{B}_n,\mbf{a},\mbf{b},\mbf{c}}^{\pm,+}
\]
in Proposition~\ref{prop:at_atrc_8v_coupling} for $\spin=\tau,\ \spin'=\tau\tau',\ K=K''=J,K'=U,\ \Lambda=\msf{B}_n$ and $\eta'=\eta^\pm$, and let $\mrm{w}_i$ be the corresponding weights in~\eqref{eq:gat_weights_1}.
Observe that $\tfrac{\mbf{a}}{\mbf{c}}=\tanh 2J$ is small enough and $\tfrac{\mbf{b}}{\mbf{c}}=\tfrac{e^{-2U}}{\cosh 2J}$ is large enough if and only if $J>0$ small enough and $-U>0$ large enough. We assume henceforth that $-U>J>0$.

By Lemma~\ref{lem:hf_law}, there exists a unique height function $h$ related to $\sigma$ by~\eqref{eq:hf_gradient} and satisfying~$h=1$ on $(\Z^2)^*\setminus\msf{B}_n^*$, and its law is given by $\hf_{\msf{B}_n,\mbf{a},\mbf{b},\mbf{c}}^{0/2,1}$. The proof is divided into three steps.

\noindent\textbf{Step 1.} For $\tfrac{\mbf{a}}{\mbf{c}}$ sufficiently small, the variance of $h$ is uniformly bounded.

Fix $N>0$. We are going to show that $h$ has uniformly bounded $N$th moment:
\[
\mbf{E}\left[\abs{h_0}^N\right]=\sum_{k\geq 1}(2k)^N\,\mbf{P}[\abs{h_0}=2k].
\]
Assume $\abs{h_0}=2k$. By the definition of height functions, this implies $\abs{h_{0'}}\geq 2k-1$, where $0'=(\tfrac{1}{2},\tfrac{1}{2})$. Since, by~\eqref{eq:hf_spin_agree}, $h$ is constant on $e\in\mbb{E}_{(\Z^2)^*}$ precisely if $\sigma^\circ$ is, we must have that $\sigma^\circ$ changes its sign at least $k-1$ times along any path in $(\Z^2)^*$ connecting $0'$ to $(\Z^2)^*\setminus \msf{B}_n^*$. However, by~\eqref{eq:at_atrc_8v_joint_law}, $\sigma^\circ$ is constant on clusters of the dual $\omega^*$, whence there exist at least $k-1$ disjoint circuits in $\omega$ within $\msf{B}_n$ that surround $0'$. One of these circuits must be of length $\ell\geq k$, of which there exist at most $\ell 4^\ell$.
Thus, by Lemma~\ref{lem:finite_energy} (as in~\eqref{eq:atrc_path_bound_1}), 
\[
\mbf{P}[\abs{h_0}=2k]\leq\sum_{\ell\geq k}\ell 4^\ell\big(1+\mrm{w}_1^{-1}\big)^{-\ell},
\]
and $\mrm{w}_1=e^{4J}-1\to 0$ as $J\to 0$. Since $\tfrac{\mbf{a}}{\mbf{c}}=\tanh 2J$, the claim is proved.

\noindent\textbf{Step 2.} Existence of $c$ and $\mcal{C}'$ when $\tfrac{\mbf{a}}{\mbf{c}}$ is sufficiently small.

Define $\mcal{C}'$ as the set of all $x'\in (\Z^2)^*$ that are connected to $(\Z^2)^*\setminus\msf{B}_n^*$ in $\omega^*$. Then $\sigma^\circ$ and thus (by~\eqref{eq:hf_spin_agree}) $h$ are constant on $\mcal{C}'$. Since $h=1$ on $(\Z^2)^*\setminus\msf{B}_n^*\subseteq\mcal{C}'$, property (ii) follows. If $\Lambda\subset (\Z^2)^*\setminus\mcal{C}'$ is connected, there must exist a circuit of length $\ell\geq\mrm{diam}(\Lambda)$ in $\omega$ within $\msf{B}_n$ that surrounds $\Lambda$. Similarly to Step 1 in the proof of Theorem~\ref{thm:antiferro}, using Lemma~\ref{lem:finite_energy}, we deduce
\[
\mbf{P}[\exists\,\Lambda\subseteq(\Z^2)^*\setminus\mcal{C}'\text{ connected with }\mrm{diam}(\Lambda)> \log n]\leq \sum_{\ell\geq\log n}\abs{\msf{B}_n}4^\ell\big(1+\mrm{w}_1^{-1}\big)^{-\ell}.
\]

\noindent\textbf{Step 3.} Existence of $\mcal{C}$ when $\tfrac{\mbf{b}}{\mbf{c}}$ is sufficiently large.

Recall the definition of the map $F_{\mrm{odd}}$ given in~\eqref{eq:sign_flip}.
Consider $(\spin,\spin')=(F_{\mrm{odd}}(\tau\tau'),\tau)$ whose law is given by $\at_{\msf{B}_n,-U,J,-J}^{+,\mrm{f}}$, and let $\widetilde{\omega}$ be distributed according to $\gat_{\msf{B}_n,-U,J,-J}^{1,\mrm{f}}$. Denote the corresponding weights in~\eqref{eq:gat_weights_1} with $K=-U,\ K'=-K''=J$ by $\widetilde{\mrm{w}}_i$.
Define $\mcal{C}$ as the set of all $x\in\Z^2$ that are connected to $\Z^2\setminus\msf{B}_n$ in $\widetilde{\omega}$. Then $F_\mrm{odd}(\tau\tau')$ is constant on $\mcal{C}$, and hence $\mbb{E}_\mcal{C}\subseteq E_{\tau\tau'}=E_{\sigma^\bullet}$. Since $h$ is constant on $e\in\mbb{E}_{\Z^2}$ if and only if $\sigma^\bullet$ is, property (i) follows. If $\Lambda\subset \Z^2\setminus\mcal{C}$ is connected, there must exist a circuit of length $\ell\geq\mrm{diam}(\Lambda)$ in $\widetilde{\omega}^*$ within $\msf{B}_n^*$ that surrounds $\Lambda$.
As in Step 2 in the proof of Theorem~\ref{thm:antiferro}, using Lemma~\ref{lem:finite_energy}, we deduce
\[
\mbf{P}[\exists\,\Lambda\subseteq\Z^2\setminus\mcal{C}\text{ connected with }\mrm{diam}(\Lambda)> \log n]\leq \sum_{\ell\geq\log n}\abs{\msf{B}_n}4^\ell\big(1+\tfrac{\widetilde{\mrm{w}}_1}{2}\big)^{-\ell},
\]
and $\widetilde{\mrm{w}}_1=e^{-2(J+U)}-1\to\infty$ as $J+U\to -\infty$. Since $\tfrac{\mbf{b}}{\mbf{c}}=\tfrac{e^{-2U}}{\cosh 2J}$, the claim follows.
\end{proof}

\section{Subcritical sharpness in the AT model}\label{sec:sharpness}
Let $d\geq 2$. In the ferromagnetic regime $J,J',U\geq 0$, the AT model enjoys monotonicity properties in the parameters. For instance, the \emph{magnetisation} at the origin $\langle\tau_0\rangle_{\msf{B}_n,J,J',U}^{+,+}$ is non-decreasing in $J,J',U$, which follows from the second Griffiths' inequality~\cite{Gri67} proved in this generality in~\cite{KelShe68} (see also~\cite{PfiVel97}). For fixed $J,J',U>0$, this implies existence of $\beta_\mrm{c}=\beta_\mrm{c}(d,J,J',U)\in [0,\infty]$ (a priori possibly trivial) such that 
\begin{equation}
\langle\tau_0\rangle_{\msf{B}_n,\beta}^{+,+}
\begin{cases}
		\xrightarrow{n\to\infty} 0 & \text{if }\beta<\beta_\mrm{c},\\
		\geq c_\beta>0 & \text{if }\beta>\beta_\mrm{c},
\end{cases}
\end{equation}
where $\langle\,\cdot\,\rangle_{\msf{B}_n,\beta}^{+,+}$ is the expectation operator with respect to the AT measure with parameters $\beta (J,J',U)$.
Then, the high-temperature expansion and a Peierls' argument~\cite{Pei36} (or comparison to the Ising model) show that $\beta_\mrm{c}\in (0,\infty)$, meaning that the model undergoes a non-trivial order-disorder phase transition.
The value of $\beta_\mrm{c}$ is known explicitly only when $d=2$ and $J=J'\geq U>0$, where it coincides with the self-dual point~\cite{AouDobGla24}.
However, the strategy in~\cite{Gri95b} (see also~\cite{HanKla23}) allows to show that the function $(J,J',U)\mapsto \beta_\mrm{c}$ is continuous.
The transition is said to be (subcritically) \emph{sharp} if,
\[
\text{for } \beta<\beta_\mrm{c}, \text{ there exists }c_\beta=c_\beta(d,J,J',U)>0 \text{ such that }\langle\tau_0\rangle_{\msf{B}_n,\beta}^{+,+}\leq e^{-c_\beta n}.
\] 
The general approach~\cite{DumRaoTas19} also applies in this setting provided a graphical representation admits the following: \emph{FKG lattice condition, maximal boundary condition, Russo-type inequality}, see Section~\ref{sec:sharp_af} for details.
When $J\geq\min\{J',U\}\geq 0$, these conditions are satisfied, see~\cite[Appendix A]{AouDobGla24} for the isotropic case. When $0\leq J<\min\{J',U\}$, sharpness for $\tau$ can be derived from the same statement for $\tau'$ and $\tau\tau'$.
Here we focus on the regime $U<0$. 

\subsection{The case when $U$ is negative: proof of Theorem~\ref{thm:sharp}}\label{sec:sharp_af}
By Lemma~\ref{lem:gat_fkg}, even when $U<0$, there is still a regime in which $\gat$ satisfies~\eqref{eq:fkg_lattice}. However, one difficulty in proving sharpness is attributed to the fact that $\tau$ and $\tau'$ are negatively correlated~\cite{PfiVel97}, e.g., for $n\geq 1$ and $x\in\msf{B}_n$,
\begin{equation}\label{eq:neg_gks}
\langle\tau_x\tau'_x\rangle_{\msf{B}_n,J,J',U}^{+,+}\leq \langle\tau_x\rangle_{\msf{B}_n,J,J',U}^{+,+}\, \langle\tau'_x\rangle_{\msf{B}_n,J,J',U}^{+,+},
\end{equation}
which is strongly related to monotonicity in the parameters in the sense of stochastic domination, and which plays a role when deriving the Russo-type inequality.
Below, we present a way to circumvent this issue by applying the strategy of~\cite{DumRaoTas19} along suitable curves in the AT phase diagram when $U<0$.
From now on, we restrict to the regime~\eqref{eq:af_fkg_regime} where $U<0$ and where the graphical representation satisfies~\eqref{eq:fkg_lattice}, which corresponds to
\begin{equation}
\big(\min\{J,J'\}>0>U \quad\text{and}\quad\tanh U\geq -\tanh J \tanh J'\,\big)\quad\Leftrightarrow\quad \mrm{w}_2<\mrm{w}_3\leq 1,
\end{equation}
where the weights $\mrm{w}_i$ are given by~\eqref{eq:gat_weights_1} with $(K,K',K'')=(J,J',U)$. Indeed, we have $U<0$ precisely if $\mrm{w}_2<\mrm{w}_3$, and the other conditions are those in~\eqref{eq:gat_fkg_regime}.
See Figure~\ref{fig:curves} for an illustration of the corresponding part of the phase diagram.

Before proceeding with the proof of Theorem~\ref{thm:sharp} via the argument in~\cite{DumRaoTas19}, it is necessary to identify the boundary condition $(1,\mrm{f})$ as maximal and to validate a certain Russo-type inequality along suitable curves that cover~\eqref{eq:af_fkg_regime}.

\paragraph{Maximal boundary condition and infinite-volume limit.}
In the regime~\eqref{eq:af_fkg_regime}, the boundary condition $(1,\mrm{f})$ is maximal for the graphical representation in the following sense.
\begin{lemma}\label{lem:max_bc}
Let $J,J',U$ satisfy~\eqref{eq:af_fkg_regime}. Let $\Lambda\subset\Delta\subset\Z^d$ be finite, and set $E_\Lambda=\overline{\mbb{E}}_\Lambda$ and $E_\Delta=\overline{\mbb{E}}_\Delta$. For any $\xi\in\Omega_{E_\Delta}^1$,
\begin{equation}\label{eq:max_bc}
\gat_{\Delta,J,J',U}^{1,\mrm{f}}[\omega_{E_\Lambda}\in\cdot \mid\omega_{E_\Delta\setminus E_\Lambda}=\xi_{E_\Delta\setminus E_\Lambda}]
\leq_\mrm{st}
\gat_{\Lambda,J,J',U}^{1,\mrm{f}},
\end{equation}
and in particular 
\begin{equation}\label{eq:mon_domain}
\gat_{\Delta,J,J',U}^{1,\mrm{f}}[\omega_{E_\Lambda}\in\cdot\,]
\leq_\mrm{st}
\gat_{\Lambda,J,J',U}^{1,\mrm{f}}.
\end{equation}
\end{lemma}
The proof is via the Holley criterion~\cite{Hol74} and is given in Appendix~\ref{sec:proof_max_bc}.
Monotonicity in the domain~\eqref{eq:mon_domain} implies existence of the weak limit 
\begin{equation}\label{eq:inf_vol_limit_af_atrc}
\gat_{J,J',U}^{1,\mrm{f}}:=\lim_{\Lambda\nearrow\Z^d}\gat_{\Lambda,J,J',U}^{1,\mrm{f}},
\end{equation}
and that it is translation-invariant and tail-trivial, and hence mixing and ergodic (see, e.g.,~\cite[Chapter 4.3]{Gri06}).

\paragraph{Russo-type inequality along suitable curves.}
We introduce a family of curves covering the regime~\eqref{eq:af_fkg_regime} along which the weights $\mrm{w}_2,\mrm{w}_3$ in~\eqref{eq:gat_weights_1} remain constant, whereas $\mrm{w}_1$ covers $(0,\infty)$. A collection of these curves is depicted in Figure~\ref{fig:curves}.

Given a curve $\gamma:(0,1)\to\R^3$ with values in~\eqref{eq:af_fkg_regime}, let $\mrm{w}_i(\gamma(\beta))$ denote the weights in~\eqref{eq:gat_weights_1} evaluated at $(K,K',K'')=\gamma(\beta)$, and write $\mrm{w}_i(\gamma)$ for the corresponding function. Moreover, we write $\gat_{\Lambda,\gamma(\beta)}^{\#,\eta'}$ for the measure in~\eqref{eq:gat_law_1} with coupling constants given by $\gamma(\beta)$. We use the same notational convention for their infinite volume limit.

\begin{lemma}\label{lem:curves}
There exists a family of disjoint smooth curves $\gamma_{\kappa,\kappa'}:(0,1)\to\R^3,\ 0<\kappa<\kappa'\leq 1,$ that satisfies the following properties: 
\begin{enumerate}[label=(\roman*)]
\item $J,J',U$ satisfy~\eqref{eq:af_fkg_regime} precisely if there exist $\kappa,\kappa'$ and $\beta$ such that $\gamma_{\kappa,\kappa'}(\beta)=(J,J',U)$,
\item for any $\kappa,\kappa'$, the weights $\mrm{w}_2(\gamma_{\kappa,\kappa'})$ and $\mrm{w}_3(\gamma_{\kappa,\kappa'})$ are constant $\kappa$ and $\kappa'$, respectively,
\item for any $\kappa,\kappa',\ \mrm{w}_1(\gamma_{\kappa,\kappa'}):(0,1)\to (0,\infty)$ is a bijection, and there exists $\varepsilon=\varepsilon(\kappa,\kappa')>0$ such that
\[
\frac{\mrm{d}}{\mrm{d}\beta}\log\mrm{w}_1(\gamma_{\kappa,\kappa'}(\beta))\geq\varepsilon.
\]
\end{enumerate}
\end{lemma}
The proof is straightforward calculus and is given in Appendix~\ref{sec:proof_curves}.
By construction, the following Russo-type inequality holds along each curve $\gamma_{\kappa,\kappa'}$. 

\begin{lemma}\label{lem:russo_ineq}
Let $\kappa,\kappa'\in (0,1]$ with $\kappa<\kappa'$. There exists $\varepsilon=\varepsilon(\kappa,\kappa')>0$ such that the following holds. For $\Lambda\subset\Z^d$ finite, $E=\overline{\mbb{E}}_\Lambda$ and any increasing random variable $X:\Omega_E^1\to [0,\infty)$, 
\begin{equation}\label{eq:russo_ineq}
\frac{\mrm{d}}{\mrm{d} \beta}\mbf{E}_{\Lambda,\gamma_{\kappa,\kappa'}(\beta)}^{1,\mrm{f}}[X]\geq\varepsilon\,\mbf{Cov}_{\Lambda,\gamma_{\kappa,\kappa'}(\beta)}^{1,\mrm{f}}(X,\abs{\omega_E}),
\end{equation}
where the expectation and covariance are taken with respect to $\gat_{\Lambda,\gamma_{\kappa,\kappa'}(\beta)}^{1,\mrm{f}}$.
In particular, $\mbf{E}_{\Lambda,\gamma_{\kappa,\kappa'}(\beta)}^{1,\mrm{f}}[X]$ is non-decreasing in $\beta$.
\end{lemma}

\begin{proof}
Fix $\kappa,\kappa',\Lambda,E$ as in the statement. For $X:\Omega_E^1\to [0,\infty)$, define
\[
Z_\beta(X):=\sum_{\omega\in\Omega_E^1}\sum_{\spin'\in\Sigma_\Lambda^{\mrm{f}}}X(\omega)\,\mrm{w}_1^{\abs{\omega_E}}\,\mrm{w}_2^{\abs{E_{\spin'}\setminus\omega}}\,\mrm{w}_3^{\abs{E_{\spin'}\cap\omega}}\,2^{k_{\Lambda}(\omega)},
\]
where $\mrm{w}_i=\mrm{w}_i(\gamma_{\kappa,\kappa'}(\beta))$. Then, since $\mrm{w}_2$ and $\mrm{w}_3$ are constant in $\beta$ by Lemma~\ref{lem:curves}(ii),
\[
\frac{\mrm{d}}{\mrm{d}\beta}Z_\beta(X)=Z_\beta(Xf_\beta)\quad\text{with}\quad f_\beta(\omega)=\big(\tfrac{\mrm{d}}{\mrm{d}\beta}\log\mrm{w}_1\big)\,\abs{\omega_E}.
\]
Differentiating $\mbf{E}_{\Lambda,\gamma_{\kappa,\kappa'}(\beta)}^{1,\mrm{f}}[X]=\tfrac{Z_\beta(X)}{Z_\beta(1)}$ and applying the quotient rule gives that the left side of~\eqref{eq:russo_ineq} coincides with $\mbf{Cov}_{\Lambda,\gamma_{\kappa,\kappa'}(\beta)}^{1,\mrm{f}}(X,f_\beta)$. 
By Lemma~\ref{lem:curves}(iii), there exists $\varepsilon=\varepsilon(\kappa,\kappa')>0$ such that $\tfrac{\mrm{d}}{\mrm{d}\beta}\log\mrm{w}_1\geq\varepsilon$.
Moreover, the covariance of $X$ and $\abs{\omega_E}$ is non-negative since both are increasing (in $\omega$) and $\gat_{\Lambda,\gamma_{\kappa,\kappa'}(\beta)}^{1,\mrm{f}}$ satisfies~\eqref{eq:fkg} by Lemma~\ref{lem:gat_fkg}. 
\end{proof}
We are now ready to run the argument of~\cite{DumRaoTas19}.
\begin{proof}[Proof of Theorem~\ref{thm:sharp}]
Consider the curves $\gamma_{\kappa,\kappa'}$ in Lemma~\ref{lem:curves}. 
Fix $d\geq 2,\ \kappa,\kappa'\in (0,1]$ with $\kappa<\kappa'$ and $\beta_0\in(0,1)$. For $k,n\geq 0$ and $\beta\in (0,\beta_0)$, define $E_k:=\overline{\mbb{E}}_{\msf{B}_k}$ and
\[
\mu_{k,\beta}:=\gat_{\msf{B}_{2k},\gamma_{\kappa,\kappa'}(\beta)}^{1,\mrm{f}},\qquad
\theta_k(\beta):=\mu_{k,\beta}[0\leftrightarrow\Z^d\setminus\msf{B}_{k-1}],\qquad
S_n(\beta):=\sum_{k=0}^{n-1}\theta_k(\beta).
\]
Then, by~\cite[Lemma 3.2]{DumRaoTas19} applied to $\mu_{n,\beta}$, 
\begin{equation}\label{eq:osss_cor}
\mbf{Cov}_{\mu_{n,\beta}}(\ind{0\leftrightarrow\Z^d\setminus\msf{B}_{n-1}},\abs{\omega_{E_{2n}}})
\geq\frac{n}{4\max\limits_{x\in\msf{B}_{n}}\sum\limits_{k=0}^{n-1}\mu_{n,\beta}[x\leftrightarrow \Z^d\setminus\msf{B}_{k-1}(x)]}\,\theta_n(\beta)(1-\theta_n(\beta)).
\end{equation}
For $x\in\msf{B}_{n}$ and $2k\leq n$, Lemma~\ref{lem:max_bc} implies ${\mu_{n,\beta}[x\leftrightarrow \Z^d\setminus\msf{B}_{k-1}(x)]}\leq \theta_k(\beta)$, and thus
\[
\sum\limits_{k=0}^{n-1}\mu_{n,\beta}[x\leftrightarrow \Z^d\setminus\msf{B}_{k-1}(x)]\leq 2 \sum_{k=0}^{\lfloor n/2 \rfloor}\mu_{n,\beta}[x\leftrightarrow \Z^d\setminus\msf{B}_{k-1}(x)]\leq 2 S_n(\beta).
\]
Plugging this in~\eqref{eq:osss_cor} and combining the result with Lemma~\ref{lem:russo_ineq}, we obtain
\begin{equation*}
\theta_n'\geq \varepsilon\,\frac{n}{8 S_n}\,\theta_n\,(1-\theta_n),
\end{equation*}
where $\varepsilon=\varepsilon(\kappa,\kappa')>0$. By inclusion of events and Lemmata~\ref{lem:max_bc} and~\ref{lem:russo_ineq}, $\theta_n(\beta)\leq\theta_0(\beta)\leq\theta_0(\beta_0)<1$. We deduce that there exists $c=c(\kappa,\kappa',d,\beta_0)>0$ such that 
\begin{equation*}
\theta_n'\geq c\frac{n}{S_n}\theta_n. 
\end{equation*}
By~\cite[Lemma 3.1]{DumRaoTas19} applied to $f_n=\theta_n/c$, there exists $\beta_1=\beta_1(d,\kappa,\kappa')\in [0,\beta_0]$ such that,
\begin{itemize}
\item for $\beta<\beta_1$, there exists $c_\beta=c_\beta(d,\kappa,\kappa')>0$ such that $\theta_n(\beta)\leq e^{-c_\beta n}$,
\item for $\beta>\beta_1$, $\theta(\beta)=\lim_{n\to\infty}\theta_n(\beta)$ satisfies $\theta(\beta)\geq c(\beta-\beta_1)$,
\end{itemize}
where $\lim_{n\to\infty}\theta_n(\beta)=\gat_{\gamma_{\kappa,\kappa'}(\beta)}^{1,\mrm{f}}[0\leftrightarrow\infty]$ as a consequence of~\eqref{eq:mon_domain} (see, e.g.,~\cite[Proposition 5.11]{Gri06}).
By Lemma~\ref{lem:curves}(iii), the weight $\mrm{w}_1(\gamma_{\kappa,\kappa'}(\beta))$ in~\eqref{eq:gat_weights_1} tends to 0 as $\beta\searrow 0$, and it tends to infinity  as $\beta\nearrow 1$. Thus, taking $\beta_0$ close enough to $1$, Lemma~\ref{lem:finite_energy} guarantees that $\beta_1\in (0,\beta_0)$. Finally, Lemma~\ref{lem:max_bc} implies
\[
\theta(\beta)\leq\gat_{\msf{B}_{2n},\gamma_{\kappa,\kappa'}(\beta)}^{1,\mrm{f}}[0\leftrightarrow\Z^d\setminus\msf{B}_{2n}]\leq \theta_n(\beta).
\]
Using the consequence, Corollary~\ref{cor:gat_corr_conn}, of the coupling in Proposition~\ref{prop:at_gat_laws}(i), we derive the statement of the theorem for $n$ even, and the case $n$ odd is treated analogously.
\end{proof}

\subsection{The isotropic case when $U$ is negative}\label{sec:sharp_af_iso}
One of the main difficulties in proving sharpness in the isotropic case when $U<0$ is that there is no known monotonicity in the parameters in the sense of stochastic domination, which is closely related to~\eqref{eq:neg_gks}. We conjecture the following property which, in conjunction with Theorem~\ref{thm:sharp}, implies sharpness (as in Theorem~\ref{thm:sharp}) along all strictly increasing curves in the isotropic phase diagram when $d=2$. Furthermore, the transition must occur at the self-dual curve defined by $\sinh 2J=e^{-2U}$~\cite{MitSte71}, as explained below.

\begin{conjecture}\label{con:mon_plus_af}
Let $\gamma:(0,1)\to \R^3$ be a curve in~\eqref{eq:af_fkg_regime} with all components strictly increasing. Then, for $\beta_1,\beta_2\in (0,1)$ with $\beta_1<\beta_2$, there exists $\varepsilon>0$ such that 
\begin{equation*}
\langle\tau_0\rangle_{\msf{B}_n,J_1,J'_1,U_1}^{+,\mrm{f}}
\leq
\langle\tau_0\rangle_{\msf{B}_n,J_2,J'_2,U_2}^{+,\mrm{f}}
\quad\text{whenever}\quad
(J_i,J'_i,U_i)\in \mcal{B}_\varepsilon(\gamma(\beta_i)),\ i=1,2,
\end{equation*}
where $\mcal{B}_r(x)$ denotes the three-dimensional Euclidean ball of radius $r>0$ and centre $x\in\R^3$. 
\end{conjecture}

We will only sketch how to use this property to derive sharpness in the isotropic case from Theorem~\ref{thm:sharp} when $d=2$. The main task is to show that, if the curve $\gamma_{\kappa,\kappa'}$ from Lemma~\ref{lem:curves} passes through the self-dual curve in the isotropic phase diagram, then the transition point $\gamma_{\kappa,\kappa'}(\beta_\mrm{c}(\kappa,\kappa'))$ from Theorem~\ref{thm:sharp} coincides with the intersection point with the self-dual curve.
Indeed, the argument in the proof of~\cite[Theorem 1.12]{Dum17a} shows that, on the curve $\gamma_{\kappa,\kappa'}$, there are only countably many points at which there is more than one infinite-volume \emph{Gibbs measure}.
It is classical (see, e.g., the proof of~\cite[Theorem 3]{AouDobGla24}) that this property, together with sharpness, Theorem~\ref{thm:sharp}, and monotonicity along $\gamma_{\kappa,\kappa'}$, Lemma~\ref{lem:russo_ineq}, implies that the self-dual point coincides with the transition point.
Finally, continuity of the curves, Conjecture~\ref{con:mon_plus_af} and Theorem~\ref{thm:sharp} easily imply the desired sharpness in the isotropic model.

We also mention that Conjecture~\ref{con:mon_plus_af} is valid in the ferromagnetic case, that is, when $\gamma$ takes values in $(0,\infty)^3$. Indeed, differentiating with respect to $J,\,J'$ and $U$ and using the second Griffiths' inequality~\cite{KelShe68}, one sees that the magnetisation at the origin is non-decreasing in each of the coupling constants. 

\paragraph{Indications of monotonicity.}
We conclude this section by presenting indications of monotonicity of the graphical representation along suitable curves that cover the isotropic part of~\eqref{eq:af_fkg_regime}, by means of the representation in Proposition~\ref{prop:atrc_pair_coupling}.
Assume henceforth that $J=J',\,U$ satisfy~\eqref{eq:af_fkg_regime}, which is equivalent to
\begin{equation}\label{eq:iso_af_fkg_regime}\tag{iso-AF-FKG}
J=J'>0>U\quad\text{and}\quad \cosh{2J}\geq e^{-2U}.
\end{equation}
When $K=K'=J,\,K''=U$ and $(\eta,\eta')=(+,\mrm{f})$, equation~\eqref{eq:atrc_pair_law} may be written as
\begin{equation}\label{eq:atrc_pair_law_iso}
\atrc_{\Lambda,J,J,U}^{1,0}[\omega,\omega']\propto \mrm{u}_1^{\abs{\omega_E}+\abs{\omega'_E}}\,\big(\tfrac{\mrm{u}_3}{\mrm{u}_1^2}\big)^{\abs{\omega_E\cap\omega'_E}}\,2^{k_\Lambda(\omega)+k_\Lambda(\omega')}\,\ind{\omega\in\Omega_E^1}\,\ind{\omega'\in\Omega_E^0},
\end{equation}
where the weights are given by~\eqref{eq:gat_weights_1} and~\eqref{eq:atrc_weights_1}, and they satisfy
\begin{equation}\label{eq:atrc_weights_2}
\mrm{u}_1=e^{2(J-U)}-1
\quad\text{and}\quad
\hat{\mrm{u}}_3:=\tfrac{\mrm{u}_3}{\mrm{u}_1^2}=\tfrac{e^{4J}-2e^{2(J-U)}+1}{(e^{2(J-U)}-1)^2}.
\end{equation}
Observe that $\hat{\mrm{u}}_3<1$ precisely if $U<0$.
Now, with respect to this measure, consider the normalised expectation
\begin{equation}\label{eq:edge_averages}
\frac{1}{\abs{E_n}}\mbf{E}_{\msf{B}_n,J,J,U}^{1,0}\left[\abs{\omega_{E_n}}+\abs{\omega'_{E_n}}\right]
=\frac{1}{\abs{E_n}}\left(\mbf{E}_{\msf{B}_n,J,J,U}^{1,\mrm{f}}\left[\abs{\omega_{E_n}}\right]+\mbf{E}_{\msf{B}_n,J,J,U}^{0,+}\left[\abs{\omega_{E_n}}\right]\right),
\end{equation}
where we chose $\Lambda=\msf{B}_n$ and $E_n=\overline{\mbb{E}}_{\msf{B}_n}$, and where the expectations on the right side are taken with respect to $\gat_{\msf{B}_n,J,J,U}^{1,\mrm{f}}$ and $\gat_{\msf{B}_n,J,J,U}^{0,+}$, respectively, see Proposition~\ref{prop:atrc_pair_coupling}. It is easy to see that the existence of the limit~\eqref{eq:inf_vol_limit_af_atrc} and the fact that it is translation-invariant imply 
\[
\frac{1}{\abs{E_n}}\mbf{E}_{\msf{B}_n,J,J,U}^{1,\mrm{f}}\left[\abs{\omega_{E_n}}\right]
\quad\xrightarrow{n\to\infty}\quad\gat_{J,J,U}^{1,\mrm{f}}[\omega_e=1],
\]
where $e$ is any fixed edge of $\Z^d$. On the other hand, a proof exactly analogous to that of~Lemma~\ref{lem:max_bc} shows that the boundary condition $(0,+)$ is minimal in the analogous sense, and hence the limit
\begin{equation}
\gat_{J,J,U}^{0,+}:=\lim_{\Lambda\nearrow\Z^d}\gat_{\Lambda,J,J,U}^{0,+}
\end{equation}
exists, and it is translation-invariant and tail-trivial. Altogether, we deduce that 
\begin{equation}\label{eq:conv_edge_averages}
\frac{1}{\abs{E_n}}\mbf{E}_{\msf{B}_n,J,J,U}^{1,0}\left[\abs{\omega_{E_n}}+\abs{\omega'_{E_n}}\right]
\quad\xrightarrow{n\to\infty}\quad
\gat_{J,J,U}^{1,\mrm{f}}[\omega_e=1]+\gat_{J,J,U}^{0,+}[\omega_e=1].
\end{equation}
We will consider the measure in~\eqref{eq:atrc_pair_law_iso} and the weights in~\eqref{eq:atrc_weights_2} with parameters given by curves taking values in~\eqref{eq:iso_af_fkg_regime}, using the same notational conventions as in the previous section.
Similarly to Lemma~\ref{lem:curves}, there exists a family of curves $\widehat{\gamma}_\kappa:(0,1)\to\R^3,\,\kappa\in (0,1),$ covering~\eqref{eq:iso_af_fkg_regime}, along which the weight $\hat{\mrm{u}}_3(\widehat{\gamma}_\kappa)$ in~\eqref{eq:atrc_weights_2} is constant $\kappa$ whereas the weight $\mrm{u}_1(\widehat{\gamma}_\kappa)$ is increasing and covers $(0,\infty)$. Moreover, along these curves, the coupling constants $J$ and $U$ are strictly increasing and decreasing, respectively. See Figure~\ref{fig:curves} for an illustration. Differentiating the expectation on the left side of~\eqref{eq:edge_averages} along these curves, we obtain
\[
\frac{\mrm{d}}{\mrm{d} \beta}\mbf{E}_{\msf{B}_n,\widehat{\gamma}_\kappa(\beta)}^{1,0}\left[\abs{\omega_{E_n}}+\abs{\omega'_{E_n}}\right]=\big(\tfrac{\mrm{d}}{\mrm{d}\beta}\log\mrm{u}_1(\widehat{\gamma}_\kappa(\beta))\big)\,\mbf{Var}_{\msf{B}_n,\widehat{\gamma}_\kappa(\beta)}^{1,0}\left(\abs{\omega_{E_n}}+\abs{\omega'_{E_n}}\right)\geq 0.
\]
Let us summarise our findings in the following proposition.
\begin{proposition}\label{prop:mon_gat_edge_densities}
There exists a family of disjoint smooth curves $\widehat{\gamma}_\kappa:(0,1)\to\R^3,\,\kappa\in (0,1),$ covering~\eqref{eq:iso_af_fkg_regime} along which the sum of the infinite-volume edge-densities 
\[
\gat_{\widehat{\gamma}_\kappa(\beta)}^{1,\mrm{f}}[\omega_e=1]+\gat_{\widehat{\gamma}_\kappa(\beta)}^{0,+}[\omega_e=1]
\]
is non-decreasing in $\beta$. For each $\kappa\in(0,1)$, $\mrm{u}_1(\widehat{\gamma}_\kappa):(0,1)\to(0,\infty)$ is strictly increasing and bijective, whereas $\hat{\mrm{u}}_3(\widehat{\gamma}_\kappa)$ is constant $\kappa$. Furthermore, along $\widehat{\gamma}_\kappa$, the coupling constants $J$ and $U$ are strictly increasing and decreasing, respectively.
\end{proposition}
Moreover, the Holley criterion~\cite{Hol74} allows to show that, for the same curves $\widehat{\gamma}_\kappa$, the measures $\gat_{\Lambda,\widehat{\gamma}_\kappa(\beta)}^{1,\mrm{f}}$  satisfy a kind of `jump monotonicity' in the stochastic sense, and it provides insight into why the regime $U<0$ is problematic.

\begin{proposition}\label{prop:jump_mon_gat}
Let $\kappa\in (0,1)$, and consider the curve $\widehat{\gamma}_\kappa$ from Proposition~\ref{prop:mon_gat_edge_densities}. Let $\beta_1,\beta_2\in (0,1)$ with $\mrm{u}_1(\widehat{\gamma}_\kappa(\beta_1)) \leq \kappa\,\mrm{u}_1(\widehat{\gamma}_\kappa(\beta_2))$. Then, for any finite $\Lambda\subset\Z^d$, 
\begin{equation}\label{eq:mon_iso_curves}
\gat_{\Lambda,\widehat{\gamma}_\kappa(\beta_1)}^{1,\mrm{f}}\leq_{\mrm{st}}\gat_{\Lambda,\widehat{\gamma}_\kappa(\beta_2)}^{1,\mrm{f}}.
\end{equation}
\end{proposition}

There exist analogues of the curves $\widehat{\gamma}_\kappa$ when $U\geq 0$, in which case the weight $\hat{\mrm{u}}_3$ in~\eqref{eq:atrc_weights_2} satisfies $\hat{\mrm{u}}_3\geq 1$. The conditions for~\eqref{eq:mon_iso_curves} imposed by the Holley criterion are then weaker and reduce to $\mrm{u}_1(\widehat{\gamma}_\kappa(\beta_1)) \leq \mrm{u}_1(\widehat{\gamma}_\kappa(\beta_2))$, that is, $\beta_1\leq\beta_2$.
In conclusion, the GAT measures are stochastically ordered along the respective curves when $U\geq 0$. However, when $U<0$, the Holley criterion guarantees stochastic ordering only if one `jumps a sufficient distance'.

The proof of Proposition~\ref{prop:jump_mon_gat} is given in Appendix~\ref{sec:jump_mon_gat}.
Notice that, together with Corollary~\ref{cor:gat_corr_conn}, the above propositions imply Proposition~\ref{prop:jump_mon_at}.

\appendix

\section{FKG lattice condition for GAT}\label{sec:proof_fkg_lattice}
\begin{proof}[Proof of Lemma~\ref{lem:gat_fkg}]
Let $K,K',K''$ with $K\geq K''$ and $K>-K''$, and omit them from the subscripts.
Fix $\Lambda\subset\Z^d$ finite, and set $E=\overline{\mbb{E}}_{\Lambda}$.
Let $\#\in\{0,1\}$ and $\eta'\in\{+,\mrm{f}\}$.

By~\cite[Theorem 2.22]{Gri06}, it suffices to show~\eqref{eq:fkg_lattice} for $\omega,\omega'\in\Omega_E^{\#}$ that agree everywhere except at two edges $e,f\in E$. Regard $\omega\in\Omega_E^{\#}$ as a set, and assume that $e,f\notin\omega$. We have to show that
\[
\gat_{\Lambda}^{\#,\eta'}[\omega\cup\{e,f\}]\,\gat_{\Lambda}^{\#,\eta'}
[\omega]\geq \gat_{\Lambda}^{\#,\eta'}[\omega\cup\{e\}]\,\gat_{\Lambda}^{\#,\eta'}[\omega\cup\{f\}].
\]
The factor $\mrm{w}_1^{\abs{\omega_E}}\,2^{k_\Lambda(\omega)}$ in~\eqref{eq:gat_law_1} satisfies the corresponding inequality (see, e.g., \cite[Theorem 3.8]{Gri06}).
Setting $Z(\zeta):=\sum_{\spin'\in\Sigma_\Lambda^{\eta'}}\mrm{w}_2^{\abs{E_{\spin'}\setminus\zeta}}\,\mrm{w}_3^{\abs{E_{\spin'}\cap\zeta}}$, it remains to show that 
\begin{equation}\label{eq:prf:fkgl_gat1}
Z(\omega\cup\{e,f\})\,Z
(\omega)\geq Z(\omega\cup\{e\})\,Z(\omega\cup\{f\}).
\end{equation}
Let $F=E\setminus\{e,f\}$, and observe that
\begin{equation}\label{eq:prf:fkgl_gat2}
Z(\zeta)=\sum_{\spin'\in\Sigma_\Lambda^{\eta'}}\mrm{w}_2^{\abs{F_{\spin'}\setminus\zeta}}\,\mrm{w}_3^{\abs{F_{\spin'}\cap\zeta}}\mrm{w}_2^{\ind{e\in E_{\spin'}\setminus\zeta}+\ind{f\in E_{\spin'}\setminus\zeta}}\mrm{w}_3^{\ind{e\in E_{\spin'}\cap\zeta}+\ind{f\in E_{\spin'}\cap\zeta}}.
\end{equation}
Define $\mu$ as the Ising measure on $\Sigma_\Lambda^{\eta'}$ given by 
\[
\mu[\spin']=\frac{1}{Z_\mu}\,\mrm{w}_2^{\abs{F_{\spin'}\setminus\omega}}\mrm{w}_3^{\abs{F_{\spin'}\cap\omega}}\qquad\text{with}\qquad Z_\mu:=\sum_{\spin'\in\Sigma_\Lambda^{\eta'}}\mrm{w}_2^{\abs{F_{\spin'}\setminus\omega}}\,\mrm{w}_3^{\abs{F_{\spin'}\cap\omega}}.
\]
Then, dividing both sides by~\(Z_\mu^2\) and inserting~\eqref{eq:prf:fkgl_gat2}, the inequality~\eqref{eq:prf:fkgl_gat1} turns into
\begin{multline*}
\mbf{E}_\mu\Big[\mrm{w}_3^{\ind{e\in E_{\spin'}}+\ind{f\in E_{\spin'}}}\Big]\,
\mbf{E}_\mu\Big[\mrm{w}_2^{\ind{e\in E_{\spin'}}+\ind{f\in E_{\spin'}}}\Big]\\
\geq
\mbf{E}_\mu\Big[\mrm{w}_2^{\ind{f\in E_{\spin'}}}\mrm{w}_3^{\ind{e\in E_{\spin'}}}\Big]\,
\mbf{E}_\mu\Big[\mrm{w}_2^{\ind{e\in E_{\spin'}}}\mrm{w}_3^{\ind{f\in E_{\spin'}}}\Big].
\end{multline*}
Applying elementary manipulations, this inequality is equivalent to
\[
\big(\mu[e,f\notin E_{\spin'}]\,\mu[e,f\in E_{\spin'}]-\mu[e\in E_{\spin'},f\notin E_{\spin'}]\,\mu[e\notin E_{\spin'},f\in E_{\spin'}]\big)\,(\mrm{w}_2-\mrm{w}_3)^2\geq 0,
\]
which is a consequence of Griffiths' second inequality~\cite{Gri67}, provided that $\mrm{w}_2,\mrm{w}_3\leq 1$. The latter is equivalent to~\eqref{eq:gat_fkg_restrictions}. Indeed, $\mrm{w}_2\leq 1$ precisely if $K'\geq K''$, and $\mrm{w}_3\leq 1$ precisely if 
\[
e^{-2K''}\leq \frac{1+e^{-2(K+K')}}{e^{-2K}+e^{-2K'}}.
\]
Applying $x\mapsto \tfrac{x-1}{x+1}$, the above inequality turns into~$-\tanh K''\leq \tanh K \tanh K'$.
\end{proof}

\section{Maximal boundary condition for GAT}\label{sec:proof_max_bc}
\begin{proof}[Proof of Lemma~\ref{lem:max_bc}]
Fix $J,J',U,\Lambda,\Delta,E_\Lambda,E_\Delta,\xi$ as in the statement, and define
\[
\mu_1:=\gat_{\Delta,J,J',U}^{1,\mrm{f}}[\omega_{E_\Lambda}\in\cdot\mid\omega_{E_\Delta\setminus E_\Lambda}=\xi_{E_\Delta\setminus E_\Lambda}]\quad\text{and}\quad\mu_2:=\gat_{\Lambda,J,J',U}^{1,\mrm{f}}.
\]
By~\cite[Theorem 2.6]{Gri06}, since $\mu_2$ satisfies~\eqref{eq:fkg_lattice} by Lemma~\ref{lem:gat_fkg}, it suffices to show that, for any $e\in E_\Lambda$ and any $\omega\in\Omega_{E_\Lambda}$ (regarded as a set) with $e\notin\omega$,
\begin{equation}\label{eq:holley1}
\frac{\mu_1[\omega\cup\{e\}]}{\mu_1[\omega]}\leq \frac{\mu_2[\omega\cup\{e\}]}{\mu_2[\omega]}.
\end{equation}
For $\zeta\in\Omega_{\Z^d}$, define
\[
Z_\Delta^{\mrm{f}}(\zeta)=\sum_{\spin'\in\Sigma_\Delta^{\mrm{f}}}\mrm{w}_2^{\abs{(E_\Delta)_{\spin'}\setminus\zeta}}\mrm{w}_3^{\abs{(E_\Delta)_{\spin'}\cap\zeta}}
\quad\text{and}\quad
Z_\Lambda^{\mrm{f}}(\zeta)=\sum_{\spin'\in\Sigma_\Lambda^{\mrm{f}}}\mrm{w}_2^{\abs{(E_\Lambda)_{\spin'}\setminus\zeta}}\mrm{w}_3^{\abs{(E_\Lambda)_{\spin'}\cap\zeta}}.
\]
After cancellations, the left side of~\eqref{eq:holley1} reduces to 
\begin{equation}\label{eq:lhs1}
\mrm{w}_1 2^{k_\Delta(\omega\cup\{e\}\cup\xi_{\mbb{E}_{\Z^d}\setminus E_\Lambda})-k_\Delta(\omega\cup\xi_{\mbb{E}_{\Z^d}\setminus E_\Lambda})}\frac{Z_\Delta^{\mrm{f}}(\omega\cup\{e\}\cup\xi_{\mbb{E}_{\Z^d}\setminus E_\Lambda})}{Z_\Delta^{\mrm{f}}(\omega\cup\xi_{\mbb{E}_{\Z^d}\setminus E_\Lambda})}.
\end{equation}
Similarly, the right side of~\eqref{eq:holley1} reduces to
\begin{equation}\label{eq:rhs1}
\mrm{w}_1 2^{k_\Lambda(\omega\cup\{e\}\cup \mbb{E}_{\Z^d}\setminus E_\Lambda)-k_\Lambda(\omega\cup\mbb{E}_{\Z^d}\setminus E_\Lambda)}\frac{Z_\Lambda^{\mrm{f}}(\omega\cup\{e\})}{Z_\Lambda^{\mrm{f}}(\omega)}.
\end{equation}
It can easily be checked that the exponent of $2$ in~\eqref{eq:lhs1} is at most that in~\eqref{eq:rhs1}. It remains to treat the ratios of partition functions. Let $F_{\Delta}=E_\Delta\setminus\{e\}$ and $F_{\Lambda}=E_\Lambda\setminus\{e\}$, and let $\nu_1$ and $\nu_2$ be the Ising measures on $\Sigma_\Delta^{\mrm{f}}$ and $\Sigma_\Lambda^{\mrm{f}}$, respectively, given by 
\[
\nu_1[\spin']\propto \mrm{w}_2^{\abs{(F_\Delta)_{\spin'}\setminus(\omega\cup\xi_{\mbb{E}_{\Z^d}\setminus E_\Lambda})}}\mrm{w}_3^{\abs{(F_\Delta)_{\spin'}\cap(\omega\cup\xi_{\mbb{E}_{\Z^d}\setminus E_\Lambda})}}
\quad\text{and}\quad
\nu_2[\spin']\propto \mrm{w}_2^{\abs{(F_{\Lambda})_{\spin'}\setminus\omega}}\mrm{w}_3^{\abs{(F_\Lambda)_{\spin'}\cap\omega}}.
\]
Then, applying elementary manipulations analogous to those in the proof of Lemma~\ref{lem:gat_fkg} above, the ratio in~\eqref{eq:lhs1} is at most that in~\eqref{eq:rhs1} if and only if
\[
\big(\nu_1[e\notin (E_\Delta)_{\spin'}]-\nu_2[e\notin (E_\Delta)_{\spin'}]\big)(\mrm{w}_3-\mrm{w_2})\geq 0. 
\]
The first factor is non-negative as a consequence of Griffiths' second inequality~\cite{Gri67}, and the second factor is non-negative precisely if $U\leq 0$.
\end{proof}

\section{Construction of the curves}
\label{sec:proof_curves}
\begin{proof}[Proof of Lemma~\ref{lem:curves}]
For $\kappa,\kappa'\in (0,1]$ with $\kappa<\kappa'$, define $\widetilde{\gamma}_{\kappa,\kappa'}:(\beta_{\kappa,\kappa'}^-,\beta_{\kappa,\kappa'}^+)\to\R^3$ by
\begin{gather*}
\beta_{\kappa,\kappa'}^-=-\tfrac{1}{4}\log (\kappa\kappa'),\quad \beta_{\kappa,\kappa'}^+=-\tfrac{1}{2}\log \kappa,\\
\widetilde{\gamma}_{\kappa,\kappa'}(\beta)=\left(
\tfrac{1}{2}\log\left(\tfrac{\kappa'-\kappa}{\kappa\kappa'e^{2\beta}-e^{-2\beta}}\right),\,
\beta,\,
\beta+\tfrac{1}{2}\log \kappa
\right).
\end{gather*}
It is straightforward to check that these curves satisfy (i) and (ii) in the statement. Moreover, the image of $(\beta_{\kappa,\kappa'}^-,\beta_{\kappa,\kappa'}^+)$ under $\mrm{w}_1(\widetilde{\gamma}_{\kappa,\kappa'})$ is $(0,\infty)$, and  
\[
-\frac{\mrm{d}}{\mrm{d}\beta}\log\mrm{w}_1(\widetilde{\gamma}_{\kappa,\kappa'}(\beta))=\frac{4\kappa(\kappa'-\kappa)e^{-4\beta}}{(\kappa\kappa'-e^{-4\beta})(e^{-4\beta}-\kappa^2)}>\frac{4\kappa}{\kappa'-\kappa}.
\]
Define $\gamma_{\kappa,\kappa'}:(0,1)\to\R^3$ by $\gamma_{\kappa,\kappa'}(\beta)=\widetilde{\gamma}_{\kappa,\kappa'}(\beta_{\kappa,\kappa'}^+-\beta(\beta_{\kappa,\kappa'}^+-\beta_{\kappa,\kappa'}^-))$.
\end{proof}

\section{Jump monotonicity for GAT}
\label{sec:jump_mon_gat}

\begin{proof}[Proof of Proposition~\ref{prop:jump_mon_gat}]
Fix $\kappa\in (0,1)$ and $\Lambda\subset\Z^d$ finite, and set $E=\overline{\mbb{E}}_\Lambda$. As $\gat_{\Lambda,\widehat{\gamma}_\kappa(\beta_i)}^{1,\mrm{f}}$ is the marginal of $\atrc_{\Lambda,\widehat{\gamma}_\kappa(\beta_i)}^{1,0}$ on the first component by Proposition~\ref{prop:atrc_pair_coupling}, the Holley criterion~\cite[Theorem 6]{Hol74} states that it suffices to show that, for all $\xi,\xi',\zeta,\zeta'\in\Omega_E$, 
\begin{equation*}
\atrc_{\Lambda,\widehat{\gamma}_\kappa(\beta_2)}^{1,0}[\xi\vee\zeta,\xi'\vee\zeta']\,
\atrc_{\Lambda,\widehat{\gamma}_\kappa(\beta_1)}^{1,0}[\xi\wedge\zeta,\xi'\wedge\zeta']\geq 
\atrc_{\Lambda,\widehat{\gamma}_\kappa(\beta_2)}^{1,0}[\xi,\xi']\,
\atrc_{\Lambda,\widehat{\gamma}_\kappa(\beta_1)}^{1,0}[\zeta,\zeta'],
\end{equation*}
where $\vee$ and $\wedge$ denote the coordinatewise maximum and minimum, respectively.
Recall that the measures can be written as in~\eqref{eq:atrc_pair_law_iso} and~\eqref{eq:atrc_weights_2}.
The factor $2^{k(\omega)+k(\omega')}$ satisfies the corresponding inequality (see, e.g.,~\cite[Theorem 3.8]{Gri06}), whence it suffices to check the above inequality for the edge-weights $\mrm{u}_1^{\abs{\omega_E}+\abs{\omega'_E}}\,\hat{\mrm{u}}_3^{\abs{\omega_E\cap\omega'_E}}$. 
Comparing respective factors for each edge $e\in E$, we obtain the following conditions:
\begin{align*}
\mrm{u}_1(\widehat{\gamma}_\kappa(\beta_2))&\geq \mrm{u}_1(\widehat{\gamma}_\kappa(\beta_1)),
	&	
\mrm{u}_1(\widehat{\gamma}_\kappa(\beta_2))^2\,\hat{\mrm{u}}_3(\widehat{\gamma}_\kappa(\beta_2))&\geq \mrm{u}_1(\widehat{\gamma}_\kappa(\beta_1))^2\,\hat{\mrm{u}}_3(\widehat{\gamma}_\kappa(\beta_1)),\\
\mrm{u}_1(\widehat{\gamma}_\kappa(\beta_2))\,\hat{\mrm{u}}_3(\widehat{\gamma}_\kappa(\beta_2))&\geq \mrm{u}_1(\widehat{\gamma}_\kappa(\beta_1)),
	&	
\mrm{u}_1(\widehat{\gamma}_\kappa(\beta_2))\,\hat{\mrm{u}}_3(\widehat{\gamma}_\kappa(\beta_2))&\geq \mrm{u}_1(\widehat{\gamma}_\kappa(\beta_1))\,\hat{\mrm{u}}_3(\widehat{\gamma}_\kappa(\beta_1)).
\end{align*}
Since $\hat{\mrm{u}}_3(\widehat{\gamma}_\kappa(\beta_1))=\hat{\mrm{u}}_3(\widehat{\gamma}_\kappa(\beta_2))=\kappa\in (0,1)$, these conditions reduce to $\mrm{u}_1(\widehat{\gamma}_\kappa(\beta_2))\,\kappa\geq\mrm{u}_1(\widehat{\gamma}_\kappa(\beta_1))$, and the proof is complete.
\end{proof}

\bibliographystyle{amsalpha}
\bibliography{biblicomplete}

\end{document}